\documentclass[preprint]{elsarticle}
\usepackage[utf8]{inputenc}
\usepackage{pdfpages}
\usepackage[hidelinks]{hyperref}
\usepackage{url}
\usepackage{amsthm}
\usepackage{thmtools}
\usepackage{amssymb}
\usepackage{amsmath}
\usepackage{float}
\usepackage{pdflscape}
\usepackage{natbib}

\usepackage{todonotes}
\usepackage{changes}
\definechangesauthor[name=MS, color=blue!50!black]{MS}
\definechangesauthor[name=KK, color=pink!50!black]{KK}
\definechangesauthor[name=JS, color=green!70!black]{JS}

\usepackage[labelformat=simple]{subcaption}

\usepackage{tikz}
\usetikzlibrary{graphs,quotes}

\usepackage[german,ruled,vlined,noend]{algorithm2e}
\SetKwProg{Fn}{Function}{}{}
\DontPrintSemicolon
\SetKwFor{For}{for}{}{}
\SetKwFor{While}{while}{}{}

\usepackage{rotating} 
\usepackage{booktabs} 
\usepackage{mathtools}

\newtheorem{defi}{Definition}
\newtheorem{theorem}[defi]{Theorem}
\newtheorem{lem}[defi]{Lemma}

\newtheorem{rem}[defi]{Remark}
\newtheorem{kor}[defi]{Corollary}
\newtheorem{example}[defi]{Example}

\newcommand{\eps}{\ensuremath{\varepsilon}}

\newcommand {\Zgeq}{\mathbb{Z}_\geq}

\newcommand {\R}{\mathbb{R}}
\DeclareMathOperator{\lexmin}{lexmin}
\DeclareMathOperator{\lexmax}{lexmax}

\DeclareMathOperator{\sort}{sort}

\newcommand{\zulM}{X} 
\newcommand{\effM}{X_{\mathrm{eff}}} 
\newcommand{\outNd}{Y_{\mathrm{nd}}} 

\newcommand {\X}{\mathcal{X}} 
\newcommand {\M}{\mathcal{M}} 
\newcommand {\E}{E} 
\newcommand {\I}{\mathcal{I}} 
\newcommand {\IS}{\mathcal{I}_S} 
\newcommand {\J}{\mathcal{I}_2} 
\newcommand{\mioc}{MIOC} 

\makeatletter
\newcommand*\bigcdot{\mathpalette\bigcdot@{1}}
\newcommand*\bigcdot@[2]{\mathbin{\vcenter{\hbox{\scalebox{#2}{$\m@th#1\bullet$}}}}}
\makeatother

\newcommand{\gPNr}{1}
\newcommand{\oPNr}{2}
\newcommand{\rPNr}{3}

\newcommand{\esa}{ESA} 

\newcommand{\lex}{\mathrm{lex}}

\newcommand {\bigo}{\mathcal{O}}
\newcommand{\preceqq}{\;\raisebox{-0.12em}{\scalebox{1.15}{\ensuremath{\mathrel{\substack{\prec\\[-.15em]=}}}}}\;}

\makeatletter
\addtotheorempostfoothook{%
   \global\everypar{{\setbox\z@\lastbox}\global\everypar{}}%
}
\makeatother

\usepackage{numcompress}\biboptions{authoryear}

\begin{document}

\title{Multi-objective Matroid Optimization with Ordinal Weights}

\author[1]{Kathrin Klamroth}
\ead{klamroth@math.uni-wuppertal.de}

\author[1]{Michael Stiglmayr}
\ead{stiglmayr@math.uni-wuppertal.de}

\author[1]{Julia Sudhoff\corref{cor1}}
\ead{sudhoff@math.uni-wuppertal.de}

\cortext[cor1]{Corresponding author}
\address[1]{University of Wuppertal, Gaußstr. 20, 42119 Wuppertal, Germany\\[3ex] Declaration of Interest: none}

\begin{abstract}
  Bi-objective optimization problems on matroids are in general intractable and their corresponding decision problems are in general NP-hard. However, if one of the objective functions is restricted to binary cost coefficients the problem becomes efficiently solvable by an exhaustive swap algorithm. Binary cost coefficients often represent two categories and are thus a special case of ordinal coefficients that are in general non-additive.
  
  In this paper we consider ordinal objective functions with more than two categories in the context of matroid optimization. 
  We introduce several problem variants that can be distinguished w.r.t.\ their respective optimization goals, analyze their interrelations, and derive a polynomial time solution method that is based on the repeated solution of matroid intersection problems.  Numerical tests on minimum spanning tree problems and on partition matroids confirm the efficiency of the approach.
\end{abstract}


\maketitle

\section{Introduction}
Matroid optimization problems, especially the minimum spanning tree problem, are well investigated even for multi-objective optimization (cf.\ \cite{Ehrgott96} for matroids and \cite{EHRGOTT97}, \cite{hamacher1994spanning}, \cite{CHEN07}, \cite{arroyo2008grasp}, \cite{davis2008hybridised} among others for spanning trees).  
In this paper, we consider multi-objective optimization problems where, in addition to one sum objective function, one or more
ordinal objective functions are to be considered. 

Ordinal coefficients occur whenever there is no numerical value to reflect the quality or cost of an element. Consider, for example, the minimum spanning tree problem which allows to find connected networks with small overall connection costs. In addition to the (non-negative integer-valued and additive) length of an edge that directly represents the cost of, e.g., a telephone cable, the construction may involve major work that affects road traffic or even public transportation systems like tramways or trains. Roads that have already a cable canal allow for a cheap and easy addition of another cable. But if it is necessary to build a new cable canal under a street, this could lead to limitations of the traffic or even affect public transportation systems. In this situation, it is useful to categorize each possible connection as ``easy to build'' (good), ``leads to little 
problems with traffic jams'' (medium), and ``leads to major problems for public transportation'' (bad), for the time of construction. It is then hard to compare, for example, a solution with two medium edges with a solution with one bad edge, since in general these categories can not be translated into monetary values. 

Ordinal weights and ordinal objective functions have been first  introduced in \cite{SCHAFER20} for shortest path problems. Motivated by applications in civil security, edges are categorized, for example, as ``secure'', ``neutral'', or ``insecure''. \cite{SCHAFER20} introduce an ordinal preorder based on ordinal weights, analyze the complexity of the problem, and suggest a polynomial time labeling algorithm for its solution. Knapsack problems with ordinal weights are analyzed in \cite{SCHAFER:knapsack}. 
They consider a general vector dominance and two lexicographic dominance concepts and suggest a dynamic programming based solution strategy and efficient greedy methods, respectively. Moreover, an outlook to multi-objective versions of ordinal problems is provided.

In this paper we extend this concept by considering multi-objective problems that combine 
one ``classical'', sum objective function with possibly several additional ordinal objectives. We focus on multi-objective optimization problems with ordinal weights on matroids and relate the multi-objective formulation to a series of single-objective optimization problems on intersections of matroids. The latter can be efficiently solved by an algorithm from \cite{Edmonds}. For the special case of bi-objective problems and only two ordinal categories, we compare this approach to the \emph{Exhaustive Swap Algorithm  (\esa)} presented in \cite{gorski21} which is even more efficient due to the special problem structure. In Section~\ref{sec:basics} we review some basic concepts of matroid theory and multi-objective optimization.

The paper is organized as follows. In Section~\ref{sec:basics} we review the basic concepts from matroid optimization and from multi-objective optimization. We particularly focus on (partial) ordering relations for ordinal objective functions in the light of state-of-the-art references. Ordinal matroid optimization problems with only one ordinal objective function are discussed in Section~\ref{sec:ordinal}. We show that ordinal matroid optimization problems can be solved by a greedy strategy, due to their special structure. Multi-objective matroid optimization problems with one sum objective and one ordinal objective are introduced in Section~\ref{sec:problem}. Their relation to matroid intersection problems is analyzed in Section~\ref{sec:matroidInter}, yielding efficient polynomial time solution strategies for all considered problem variants. The algorithms are numerically tested and compared at randomly generated instances of graphic matroids and of partition matroids in Section~\ref{sec:numeric}, and the paper is concluded with a short outlook on future research topics in Section~\ref{sec:concl}.

\section{Preliminaries}
\label{sec:basics}
Since we combine matroid theory and multi-objective combinatorial optimization in this paper, this preliminaries section is divided into two parts.
In the first subsection we summarize some basic definitions and results of matroid theory (for a self-contained introduction to matroid theory see, for example, \cite{oxley11matroid,schrijver02combinatorial,edmonds71matroids,Schrijver13}). In the second subsection we review basic concepts of multi-objective optimization with a particular focus on dominance relations. For a general introduction into multi-objective optimization see, e.g., \cite{Ehrg05}. A survey on multi-objective combinatorial optimization is given in \cite{ehrgott00survey}.

\subsection{Matroid Theory}
Let $\E$ be a finite set and let $\I \subset 2^\E$ be a subset of the power-set of $\E$. The tuple $(\E,\I)$
is called a \emph{matroid} if and only if the following three properties hold:
\begin{align}
	&\emptyset \in \I \label{eq:memptyset}\\
	&I\in\I\text{ and } J\subseteq I \implies  J\in\I \label{eq:minclusion}\\
	&I,J\in\I\colon \vert I\vert<\vert J\vert \implies \exists j\in J\setminus I \colon I\cup\{j\}\in\I. \label{eq:mexchange}
\end{align}
The subsets $I\in \I$ are called \emph{independent sets} while the subsets $D\in2^\E\setminus \I$ are called \emph{dependent sets}. Moreover, if the tuple $(\E,\I)$ satisfies at least the conditions (1) and (2), then it is called an \emph{independence system}.

Furthermore, all inclusion-wise maximal independent sets are called \emph{bases}, and all inclusion-wise minimal dependent sets are called \emph{circuits}. In the following, we write $\X\coloneqq\{B\in\I:\nexists I\in\I:I\supsetneq B\}$ for the set of all bases of a matroid.
All bases of a matroid have the same cardinality which is referred to as the \emph{rank} of the matroid, see, for example, \cite{oxley11matroid}. 
An important characteristic of matroids is the so-called \emph{basis exchange property}: 
	\begin{align}\label{basisexprop}
		& \forall B_1,B_2\!\in\X\;\; \forall b_1\!\in\! B_1\setminus B_2\;\; \exists \,b_2\!\in\! B_2\setminus B_1\;\colon\;
		(B_1\cup\{b_2\})\setminus\{b_1\} \in \X .
	\end{align}
The \emph{restriction} of a matroid $\M=(\E,\I)$ to a subset $S\subseteq \E$ is defined as $\M-S\coloneqq(\E\setminus S,\IS)$ with $\IS=\{I\in\I:I\subseteq \E\setminus S\}$. Another way to manipulate a given matroid $\M$ is the \emph{contraction} of $\M$ by an independent set $I\in\I$ that is defined as $\M/I\coloneqq(\E\setminus I,\mathcal{K})$ with $\mathcal{K}=\{K\subseteq \E\setminus I:K\cup I\in\I\}$. 

To simplify the notation we define the set operations $S+e\coloneqq S\cup\{e\}$ and $S-f\coloneqq S\setminus\{f\}$ for $S\subseteq E$ and  $e,f\in \E$. Furthermore, let $S^c\coloneqq \E\setminus S$ denote the complement of $S$ in $\E$.

For an extensive list of examples of matroids we refer to \cite{oxley11matroid}. The following three matroids are frequently considered and will be used for illustrations and numerical tests in this paper.
	\begin{itemize}
		\item \textbf{graphic matroid}:
		Let $G=(V,E)$ be an undirected graph, then $\M=(E,\I)$  with $\I\coloneqq \{I\subseteq E\colon (V,I)\text{ contains no cycle}\}$ is a matroid.
		\item \textbf{uniform matroid}: 
		Let $\E$ be a finite set and $k\in\mathbb{N}_0$, then $\M=(\E,\I)$ with $\I\coloneqq \{I\subseteq \E\colon \vert I\vert\leq k\}$ is a matroid.
		\item \textbf{partition matroid}: 
		Let $\E=E_1\cup E_2\cup\ldots\cup E_k$ be the disjoint union of $k$ finite sets and let $u_1,\dots,u_k\geq0$ be non-negative integers. Then $\M=(\E,\I)$ with $\I\coloneqq \{I\subseteq \E\colon \vert I\cap E_i\vert\leq u_i \;\forall 1\leq i\leq k\}$ is a matroid.
	\end{itemize}

An important concept that will prove very useful in this paper is the intersection of two matroids. Consider two matroids $\M_1=(\E,\I_1)$ and $\M_2=(\E,\I_2)$ over the same ground set $\E$. Then the \emph{matroid intersection} of $\M_1$ and $\M_2$ is defined as the independence system $\M_1\cap \M_2\coloneqq (\E,\I_1 \cap\I_2)$. 

\begin{figure}[h]
        \hspace*{\fill}
        \subcaptionbox{\label{fig:mi_1}}{


\begin{tikzpicture}[scale=0.7, every node/.style={scale=0.8}]
\node[draw,circle,thick,fill=black!10] (1) at (0,2)[] {$1$};
\node[draw,circle,thick,fill=black!10] (2) at (0,0) {$2$};
\node[draw,circle,thick,fill=black!10] (3) at (1,1) {$3$};
\node[draw,circle,thick,fill=black!10] (4) at (3,1) {$4$};
\node[draw,circle,thick,fill=black!10] (5) at (4,2) {$5$};
\node[draw,circle,thick,fill=black!10] (6) at (4,0) {$6$};


\graph {
(1) --[thick,draw=black!40!red] (3);
(2) --[thick,swap, draw=black!40!red] (3);
(1) --[dotted,thick,swap,draw=black!50!green] (2);
(3) --[thick,swap,draw=black!40!red] (4);
(4) --[dotted,draw=black!50!green,thick] (5);
(4) --[dotted,thick,swap,draw=black!50!green] (6);
(5) --[thick,draw=black!40!red] (6);

};
\end{tikzpicture}

} 
        \hspace*{\fill}\hspace*{\fill}
        \subcaptionbox{\label{fig:mi_2}}{


\begin{tikzpicture}[scale=0.7, every node/.style={scale=0.8}]
\node[draw,circle,thick,fill=black!10] (1) at (0,2)[] {$1$};
\node[draw,circle,thick,fill=black!10] (2) at (0,0) {$2$};
\node[draw,circle,thick,fill=black!10] (3) at (1,1) {$3$};
\node[draw,circle,thick,fill=black!10] (4) at (3,1) {$4$};
\node[draw,circle,thick,fill=black!10] (5) at (4,2) {$5$};
\node[draw,circle,thick,fill=black!10] (6) at (4,0) {$6$};


\graph {
(1) --[thick,draw=black!40!red] (3);
(2) --[thick,swap, draw=black!40!red] (3);
(4) --[dotted,draw=black!50!green,thick] (5);
(4) --[dotted,thick,swap,draw=black!50!green] (6);
};
\end{tikzpicture}

}         
        \hspace*{\fill}\hspace*{\fill}
        \subcaptionbox{\label{fig:mi_3}}{


\begin{tikzpicture}[scale=0.7, every node/.style={scale=0.8}]
\node[draw,circle,thick,fill=black!10] (1) at (0,2)[] {$1$};
\node[draw,circle,thick,fill=black!10] (2) at (0,0) {$2$};
\node[draw,circle,thick,fill=black!10] (3) at (1,1) {$3$};
\node[draw,circle,thick,fill=black!10] (4) at (3,1) {$4$};
\node[draw,circle,thick,fill=black!10] (5) at (4,2) {$5$};
\node[draw,circle,thick,fill=black!10] (6) at (4,0) {$6$};


\graph {
(1) --[thick,draw=black!40!red] (3);
(1) --[dotted,thick,swap,draw=black!50!green] (2);
(3) --[thick,swap,draw=black!40!red] (4);
(4) --[dotted,draw=black!50!green,thick] (5);
(4) --[dotted,thick,swap,draw=black!50!green] (6);
};
\end{tikzpicture}

}  
        \hspace*{\fill}
\caption{The edges $E=E_1\cup E_2$ of the graph \(G=(V,E)\) in \subref{fig:mi_1} are the ground-set of a graphic matroid $\M_1$ and of a partition matroid $\M_2$ (where $E_1$ is the set of all green-dotted edges and $E_2$ is the set of all red-solid edges). The independent set of edges \(I\) illustrated in \subref{fig:mi_2} is
inclusion-wise maximal for $\M_1\cap \M_2$, 
however, the alternative independent set $J$ shown in \subref{fig:mi_3} has larger cardinality.}
\label{fig:noMatroid}
\end{figure}
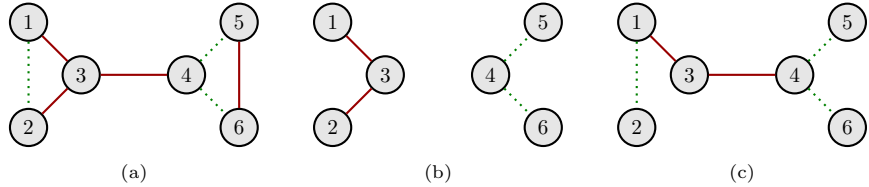

It is important to note that a matroid intersection is not necessarily a matroid itself. 
A counter example is given in Figure~\ref{fig:noMatroid}. Let \(G=(V,E_1\cup E_2)\) be a graph with an edge set $\E=E_1\cup E_2$ that is partitioned into two subsets $E_1,E_2$ (e.g., green and red edges,  respectively). Moreover, let $\M_1=(\E,\I_1)$ be the graphic matroid on \(G\) and
let $\M_2=(\E,\I_2)$ be a partition matroid with $E_1\coloneqq\{[1,2],[4,5],[4,6]\}$ the set of all green edges, $E_2\coloneqq\{[1,3],[2,3],[3,4],[5,6]\}$ the set of all red edges and $\I_2\coloneqq \{I\subseteq \E\colon \vert I \cap E_1\vert\leq 3,\, \vert I\cap E_2\vert\leq 2\}$. 
Then the set $I\coloneqq\{[1,3],[2,3],[4,5],[4,6]\}$ is an inclusion-wise maximal independent set of $\M_1\cap \M_2$ since every additional edge \(e\in \E\setminus I\) makes \(I\cup\{e\}\) dependent w.r.t.\ either \(\M_1\) or \(\M_2\). 
However, the set $J\coloneqq\{[1,2],[1,3],[3,4],[4,5],[4,6]\}$ is also a maximal independent set of $\M_1\cap \M_2$ that has larger cardinality, i.e., $\vert J\vert>\vert I\vert$. Thus, $\M_1\cap \M_2$ is not a matroid because all maximal independent sets of a matroid must have the same cardinality.
However, if at least one of the intersected matroids is a uniform matroid then the intersection is again a matroid, see \cite{oxley11matroid}.

\subsection{Multi-objective Optimization}\label{subsec:multiobjective}

In this article we consider three different concepts of efficiency for multi-objective optimization. The first two are used for general, usually sum objective functions, while the last one is particularly defined for ordinal objective functions. Throughout this section we consider minimization problems with a discrete feasible set $\zulM \subset 2^\E$ that is defined over a finite ground set $\E$.  However, all of the presented concepts can be defined analogously for maximization problems.

\paragraph{Pareto Optimality}
The Pareto concept of optimality (see again, e.g., \cite{Ehrg05}) is based on the component-wise order. Let  $y^1,y^2\in\mathbb{R}^p$ and define 
	\begin{align*}
		y^1\leqq y^2 &\,:\Longleftrightarrow\, y_i^1 \leq y_i^2,\; \quad i=1,\ldots ,p,\\
		y^1\leqslant y^2 &\,:\Longleftrightarrow\, y_i^1\leq y_i^2,\; \quad i=1,\ldots,p \;\;\text{and}\;\; y^1\neq y^2,\\
		y^1<y^2 &\,:\Longleftrightarrow\, y_i^1<y_i^2, \;\quad i=1,\ldots,p.
	\end{align*}
We say that \(y^1\) \emph{dominates} \(y^2\) whenever \(y^1\leqslant y^2\). Note that the binary relation $\leqq$ is reflexive, transitive and antisymmetric and hence a partial order. The binary relation $\leqslant$ is irreflexive, transitive and asymmetric and hence defines a strict partial order. Similarly, $<$ defines a strict partial order.

This concept of dominance is typically used when comparing outcome vectors of a \emph{multi-objective problem} \eqref{eq:MOP} that aims at minimizing a vector-valued objective function $w:\zulM\to\mathbb{R}^p$: 
\begin{equation}\label{eq:MOP}\tag{MOP}
  \begin{array}{rl}
    \min & w(x)\\
    \text{s.\,t.} & x\in \zulM . 
  \end{array}
\end{equation}
A feasible solution $\bar{x}\in \zulM$ is called \emph{efficient} or \emph{Pareto-optimal} for \eqref{eq:MOP} if there exists no $x\in \zulM$ with $w(x)\leqslant w(\bar{x})$. Moreover, a feasible solution $\hat{x}\in \zulM$ is called \emph{weakly efficient} or \emph{weakly Pareto-optimal} for \eqref{eq:MOP} if there exists no $x\in \zulM$ with $w(x)< w(\hat{x})$.
The set of all efficient solutions of a problem \eqref{eq:MOP} is called \emph{efficient set} and denoted by $\effM$. 
The image $w(\bar{x})$ of an efficient solution $\bar{x}$ is called \emph{non-dominated} outcome vector or non-dominated point. 
If $w(x')\leqslant w(x'')$ for two feasible solutions $x',x''\in \zulM$, then we say that $x'$ \emph{dominates} $x''$, in accordance with the above concept of dominance in the objective space.
The image of the efficient set is called \emph{non-dominated set} and denoted by  $\outNd\coloneqq w(\effM)$.

\paragraph{Lexicographic Optimality} The concept of lexicographic optimality assumes a specific ordering among the components of the given (outcome) vectors, i.e., the first component is more important than the second, and so on. We refer again to \cite{Ehrg05} for a more general introduction.
Let $y^1,y^2\in\mathbb{R}^p$. Then \(y^1<_{\lex} y^2\), i.\,e., \(y^1\) \emph{ lexicographically dominates} \(y^2\), 
if there is an index \(k\in\{1,\ldots,p\}\) such that \(y^1_k < y^2_k\) and \(y^1_i=y^2_i\) for all \(i\in\{1,\ldots,k-1\}\). Furthermore, we write \(y^1\leqq_{\lex} y^2\) if \(y^1<_{\lex} y^2\) or \(y^1=y^2\). Note that the lexicographic order $\leqq_{\lex}$ is a total order, i.\,e., it is reflexive, transitive and antisymmetric, and for all $y^1,y^2\in\mathbb{R}^p$ with $y^1\neq y^2$ either \(y^1 \leqq_{\lex} y^2\) or \(y^2 \leqq_{\lex} y^1\) holds. 
Consequently, we call a solution \(\bar{x}\in \zulM\) of \eqref{eq:MOP}  \emph{lexicographically optimal} if 
$w(\bar{x})\leqq_{\lex} w(x)$ for all \(x\in \zulM\). To distinguish lexicographic optimization from multi-objective optimization in the Pareto sense we write \(\lexmin\) (and \(\lexmax\) in the case of maximization problems, respectively). 

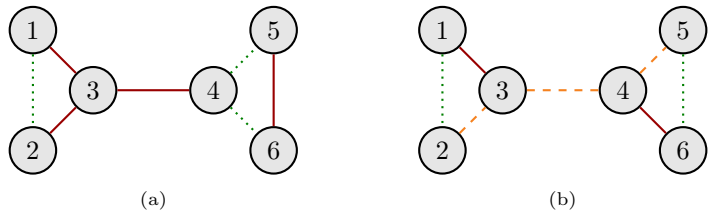
\begin{figure}
	\hspace*{\fill}
	\subcaptionbox{\label{fig:bbmpVsMop:1}}{


\begin{tikzpicture}[scale=0.8]
\node[draw,circle,thick,fill=black!10] (1) at (0,2)[] {$1$};
\node[draw,circle,thick,fill=black!10] (2) at (0,0) {$2$};
\node[draw,circle,thick,fill=black!10] (3) at (1,1) {$3$};
\node[draw,circle,thick,fill=black!10] (4) at (3,1) {$4$};
\node[draw,circle,thick,fill=black!10] (5) at (4,2) {$5$};
\node[draw,circle,thick,fill=black!10] (6) at (4,0) {$6$};


\graph {
(1) --[thick,draw=black!40!red] (3);
(2) --[thick,swap, draw=black!40!red] (3);
(1) --[dotted,thick,swap,draw=black!50!green] (2);
(3) --[thick,swap,draw=black!40!red] (4);
(4) --[dotted,draw=black!50!green,thick] (5);
(4) --[dotted,thick,swap,draw=black!50!green] (6);
(5) --[thick,draw=black!40!red] (6);

};
\end{tikzpicture}

} 
	\hspace*{\fill}
	\subcaptionbox{\label{fig:bbmpVsMop:2}}{


\begin{tikzpicture}[scale=0.8]
\node[draw,circle,thick,fill=black!10] (1) at (0,2)[] {$1$};
\node[draw,circle,thick,fill=black!10] (2) at (0,0) {$2$};
\node[draw,circle,thick,fill=black!10] (3) at (1,1) {$3$};
\node[draw,circle,thick,fill=black!10] (4) at (3,1) {$4$};
\node[draw,circle,thick,fill=black!10] (5) at (4,2) {$5$};
\node[draw,circle,thick,fill=black!10] (6) at (4,0) {$6$};


\graph {
(1) --[thick,draw=black!40!red] (3);
(2) --[dashed,thick,swap, draw=yellow!40!red] (3);
(1) --[dotted,thick,swap,draw=black!50!green] (2);
(3) --[dashed,thick,swap,draw=yellow!40!red] (4);
(4) --[dashed,draw=yellow!40!red,thick] (5);
(4) --[thick,swap,draw=black!40!red] (6);
(5) --[dotted,thick,draw=black!50!green] (6);

};
\end{tikzpicture}

}
	\hspace*{\fill}
	\caption{The edges of the graph in \subref{fig:bbmpVsMop:1} are categorized w.r.t.\ two categories(green-dotted and red-solid), c.\,f.\ Figure~\ref{fig:mi_1} above, while the edges of the same graph are categorized w.r.t.\ three categories in \subref{fig:bbmpVsMop:2} (green-dotted, orange-dashed, and red-solid).} 
	\label{fig:bbmpVsMop}	
\end{figure}

\paragraph{Ordinal Optimality} Now assume that all solutions in the feasible set $X$ have the same cardinality, i.\,e.\ \(|x|=r\) for all $x\in X$. This is, for example, satisfied in the case of a matroid optimization problem. Then   
an \emph{ordinal weight optimization problem} \eqref{eq:OWOP} can be formulated as
\begin{equation}\label{eq:OWOP}\tag{OWOP}
  \begin{array}{rl}
    \min & o(x)\\
    \text{s.\,t.} & x\in \zulM.
  \end{array}
\end{equation} 
Intuitively, the ordinal objective function $o$ assigns one out of $K$ ordered categories to each element of the ground set $E$, and hence the ordinal objective of $x$ is given by an $r$-dimensional ordinal vector. For example, in the case $K=3$ we may think of good (green), medium (orange) and bad (red) elements, where we prefer good over medium and medium over bad.
Figure~\ref{fig:bbmpVsMop} shows two examples of the ground set $E$ of a graphic matroid. While the edges in Figure~\ref{fig:bbmpVsMop:1} are assigned to only two categories (where green-dotted is better than red-solid), Figure~\ref{fig:bbmpVsMop:2} shows an example with three categories (where green-dotted is better than orange-dashed which is again better than red-solid).

Throughout this paper we assume that the components of the ordinal vectors of feasible solutions are sorted in non-decreasing order w.r.t.\ the quality of the respective categories, see Figure~\ref{fig:allBases} for an illustration. This sorting will be useful when comparing different solutions in the following. 

More formally, let $\mathcal{C}=\{\eta_1,\ldots,\eta_K\}$ be an ordinal space consisting of $K$ ordered categories, and let $o:E\to \mathcal{C}$ assign one ordinal category to each element of the ground set $E$. Moreover, (by slightly abusing the notation) let 
$o:\zulM\to \mathcal{C}^r$ be a function mapping each feasible solution to an \(r\)-dimensional ordinal vector. We assume that category \(\eta_i\) with \(i\in\{1,\ldots,K-1\}\) is strictly preferred over all categories \(\eta_{j}\) with $i<j$, which is denoted by $\eta_i \prec \eta_{j}$. Similarly, we write $\eta_i\preceq \eta_j$ whenever $i\leq j$. Moreover, the components of the objective vector $o(x)$ of a feasible solution $x\in X$ are sorted in non-decreasing order, which is denoted by $o(x)\coloneqq \sort(o(x_1),\dots,o(x_r))$. If we specify the vector $o$ for an explicit example we often write $i$ instead of $\eta_i$ for better readability.

In order to define meaningful optimality concepts for problem \eqref{eq:OWOP}, we need to compare ordinal vectors in $\mathcal{C}^r$. The following definition is based on the concept first introduced in \cite{SCHAFER20} and \cite{SCHAFER:knapsack}. Let $y^1,y^2\in\mathcal{C}^r$ be two ordinal vectors. Then we write 
\begin{align*}
	y^1\preceqq y^2 &:\Longleftrightarrow y_i^1 \preceq y_i^2,\; i=1,\ldots ,r,\\
	y^1\preccurlyeq y^2 &:\Longleftrightarrow y_i^1\preceq y_i^2,\; i=1,\ldots,r\text{ and } y^1\neq y^2,\\
	y^1\prec y^2 &:\Longleftrightarrow y_i^1\prec y_i^2, \; i=1,\ldots,r.
\end{align*}
When $y^1=o(x^1)$ and $y^2=o(x^2)$ are outcome vectors of problem \eqref{eq:OWOP}, then their components are sorted in non-decreasing order. We use the same notation in this case, and
we say that $y^1$ \emph{ordinally dominates} $y^2$ whenever $y^1\preccurlyeq y^2$. Note that what we consider here is a special case of the concept of ordinal dominance introduced in \cite{SCHAFER20} who considered the more general case when feasible solutions --  and hence their outcome vectors -- may differ w.r.t.\ their number of elements. \cite{SCHAFER20} showed that in this more general case, the binary relation $\preceqq$ is a partial preorder on the set of sorted outcome vectors of an ordinal optimization problem, i.\,e., it is reflexive and transitive. In the special case when all feasible solutions have the same number of elements, as considered in this paper, the binary relation $\preceqq$ is also antisymmetric and thus a partial order. See also \cite{SCHAFER:knapsack} for yet another perspective on ordinal efficiency. Moreover, the binary relation $\preccurlyeq$ is a strict partial order in our case, i.\,e., it is irreflexive, transitive and asymmetric. 
Thus, the concepts of (weak) ordinal efficiency and (weak) ordinal dominance can be defined in a similar way as for the case of Pareto optimality by replacing $\leqslant$ with $\preccurlyeq$ and $<$ with $\prec$, respectively.

\paragraph{Combined Orderings} 
If an optimization problem has an objective function that maps feasible solutions to outcome vectors on which several of the above orderings are combined -- this may be the case when, for example, the first $p$ components of an outcome vector represent sum objective functions that are ordered w.r.t.\ Pareto dominance, while the following $r$ objective values represent ordinal values to which ordinal dominance is applied -- then we say that a solution \(x'\) \emph{dominates} a solution \(x''\), if all objective values of $x'$ are ``at least as good'' w.r.t.\ all components with at least one strict inequality in the respective ordering concept.

\section{Ordinal Matroid Optimization}\label{sec:ordinal}
As a first step towards multi-objective ordinal optimization problems we investigate matroid optimization problems with only one ordinal objective function and show that such problems can be solved using a greedy algorithm. In slight abuse of the standard notation, we will refer to the resulting problems as ``single-objective optimization problems'', even though their objective functions are vector-valued. Similar results were obtained by \cite{SCHAFER20} and \cite{SCHAFER:knapsack} for shortest path and knapsack problems, respectively, however, for the case that a lexicographic optimization is employed on the ordinal outcome vectors. We show in the following that in the case of matroids, ordinal optimality actually coincides with lexicographic optimality, and hence a greedy algorithm always yields the ordinally non-dominated set in this case. 

\subsection{Ordinal and Lexicographic Optimality and their Interrelation}\label{subsec:singleobjmodels}

Let a matroid $\M_1=(\E,\I_1)$ with rank $r$ be given and denote the set of its bases by \(\mathcal{X}_1\).
As a first special case, consider the situation of only two categories, i.e., $K=2$. Wlog we set $\eta_1=0$ (green) and $\eta_2=1$ (red). It is easy to see that in this case problem \eqref{eq:OWOP} is equivalent to a matroid optimization problem with a ``classical'' sum objective function with binary coefficients $b:\E\to\{0,1\}$ where the cost of a basis \(B\) is the aggregated cost of all of its elements, i.\,e.\ $b(B)\coloneqq\sum_{e\in B} b(e)$. Indeed, a basis $B_1\in\mathcal{X}_1$ ordinally dominates a basis $B_2\in\mathcal{X}_1$ whenever the number of one-entries in $o(B_1)$ (i.e., red elements in $B_1$)  is smaller than that in $o(B_2)$.
This leads to a \emph{matroid problem with a binary objective function} \eqref{p:bmp}
\begin{equation}\label{p:bmp}\tag{BMP}
  \begin{array}{rl}
    \min & b(B)\\
    \text{s.\,t.} & B\in \X_1 
  \end{array}
\end{equation}
as a particularly simple special case of problem \eqref{eq:OWOP}.
When $K>2$, i.e., when more than two ordinal categories have to be considered, a simple aggregation of all categories into one single aggregated objective value is no longer meaningful. However, we will discuss two related optimization problems that are based on partial aggregation in the following. Towards this end, let problem \eqref{eq:OWOP} for the special case of matroid optimization be defined as a \emph{matroid problem with ordinal costs}, given by
\begin{equation}\label{p:mpo}\tag{MPO}
  \begin{array}{rl}
    \min & o(B)\\
    \text{s.\,t.} & B\in \X_1 .
  \end{array}
\end{equation} 
Now consider a feasible basis $B\in \X_1$. Then the information contained in the objective vector $o(B)\in\mathcal{C}^r$ can equivalently be stored in an aggregated vector $c(B)\in\Zgeq^K$ with components $c_j(B)\coloneqq\vert\{e\in B \colon o(e)=\eta_j\}\vert$ for $j=1,\dots,K$ that count the number of elements in each category in $B$. Indeed, there is a simple one-to-one correspondence between $o(B)$ and $c(B)$. We will refer to $c$ as a \emph{counting objective function} in the following. This representation is often advantageous since in general \(K\), i.\,e., the number of categories, is constant and much smaller than the dimension \(r\), i.\,e., the number of elements in a basis.
Note that since all bases have the same number of elements we have that $\sum_{j=1}^K c_j(B)=r$, and hence all outcome vectors $c(B)$, $B\in\X_1$, lie on the same hyperplane in $\R^K$. Moreover, one of the components of $c$ can be omitted without loosing any information.
 
This reformulation suggests two related lexicographic optimization problems: On the one hand, we may aim at lexicographically
maximizing the number of elements in the ``good'' categories, 
and on the other hand, we may want to lexicographically minimize the number of elements in the ``bad'' categories. 
In order to clearly distinguish between these two optimization goals, we introduce two separate variants of the counting objective $c$ denoted as $c^{\max}$ and $c^{\min}$, respectively.

\paragraph{Maximizing the Number of Good Elements}
When aiming at the maximization of the number of elements in \emph{good} categories, we can apply a lexicographic \emph{maximization} to the counting objective $c$. Thus, in this case we set
$c^{\max}(B)\coloneqq c(B)$ for $B\in\X_1$ with $c^{\max}_j(B)\coloneqq\vert\{e\in B \colon o(e)=\eta_j\}\vert$ for $j=1,\dots,K$ as defined above, and formulate problem \eqref{p:mpcmax} as
 \begin{equation}\label{p:mpcmax}\tag{MPCmax}
	\begin{array}{ll}
		\lexmax &c^{\max}(B)\\
		\text{s.\,t.}&B\in \X_1 .
	\end{array}
\end{equation}

\paragraph{Minimizing the Number of Bad Elements}
In order to lexicographically minimize the number of elements in \emph{bad} categories, we first have to bring the corresponding entries of the counting objective $c$ that represent the bad categories into the leading positions (which are always considered first in lexicographic optimization). We hence define $c^{\min}_j(B)\coloneqq c_{K-j+1}(B)$ for $j=1,\dots,K$ and for $B\in\X_1$, i.\,e., $c^{\min}_j(B)\coloneqq\vert\{e\in B \colon o(e)=\eta_{K-j+1}\}\vert$, and consider problem \eqref{p:mpc} given by
\begin{equation}\label{p:mpc}\tag{MPCmin}
	\begin{array}{ll}
		\lexmin &c^{\min}(B)\\
		\text{s.\,t.}&B\in \X_1 .
	\end{array}
\end{equation}

Figure~\ref{fig:allBases} shows an example of a graphic matroid with all of its feasible bases and their respective objective vectors $o$, $c^{\max}$ and $c^{\min}$, see also Example~\ref{ex:allBases} below.

\subsection{Interrelation Between \eqref{p:mpo}, \eqref{p:mpc} and \eqref{p:mpcmax}}\label{sec:inerrelateo}
In general, the ordinally non-dominated set of problem \eqref{p:mpo} is different from the sets of lexicographically optimal outcome vectors of the associated formulations \eqref{p:mpc} and \eqref{p:mpcmax}, respectively. This can be seen, for example, at the cases of ordinal shortest path problems (see \cite{SCHAFER20}) and ordinal knapsack problems (see \cite{SCHAFER:knapsack}). In the special case of matroids, however, these three concepts are closely related and their respective efficient and non-dominated sets coincide.

\begin{theorem}\label{thm:ParetoVSlex} 
	Let $\M_1=(\E,\I_1)$ be a matroid, let $\mathcal{X}_1$ denote the set of bases of $\M_1$, and let the functions $o$, $c^{\min}$ and $c^{\max}$ be given and defined as above. Moreover, let $B_1, B_2 \in \X_1$ be two bases of $\M_1$. Then 
	$$ \left(\vphantom{c^{\min}} o(B_1) \preccurlyeq o(B_2)\right) \Rightarrow \left(c^{\min}(B_1)<_{\lex}c^{\min}(B_2) \text{~and~} c^{\max}(B_1)>_{\lex}c^{\max}(B_2)\right),
	$$
	i.\,e., if $o(B_1)$ ordinally dominates $o(B_2)$, then $c^{\min}(B_1)$ lexicographically dominates $c^{\min}(B_2)$ and $c^{\max}(B_1)$ lexicographically dominates $c^{\max}(B_2)$.
\end{theorem}
\begin{proof}
	We prove the result for $c^{\min}$. The corresponding result for $c^{\max}$ follows analogously, noting that \eqref{p:mpc} involves lexicographic minimization while \eqref{p:mpcmax} involves lexicographic maximization.
	
	Now let $o(B_1) \preccurlyeq o(B_2)$ and assume that $c^{\min}(B_1)$ does not lexicographically dominate $c^{\min}(B_2)$. First note that $o(B_1) \preccurlyeq o(B_2)$ implies $o(B_1) \neq o(B_2)$ and hence $c^{\min}(B_1)\neq c^{\min}(B_2)$. Let $\tau\coloneqq \min\{i:c_i^{\min}(B_1)\neq c_i^{\min}(B_2)\}$ be the smallest index where $c^{\min}(B_1)$ and $c^{\min}(B_2)$ differ. Since we assumed that $c^{\min}(B_1)$ does not lexicographically dominate $c^{\min}(B_2)$, it follows that $c_{\tau}^{\min}(B_1) > c_{\tau}^{\min}(B_2)$. Thus, the vectors $o(B_1)$ and $o(B_2)$ are equal in the last $\ell\coloneqq\sum_{i=1}^{\tau-1}c_i^{\min}(B_1)=\sum_{i=1}^{\tau-1}c_i^{\min}(B_2)$ components, i.\,e., \(o_j(B_1)=o_j(B_2)\) for all \(j=K-\ell+1,\ldots, K\). Furthermore, it holds that $o_{K-\ell}(B_1) \succ o_{K-\ell}(B_2)$, which contradicts the assumption that $o(B_1)$ ordinally dominates $o(B_2)$.
\end{proof}

Note that while the proof of Theorem~\ref{thm:ParetoVSlex} relies on the fact that all feasible solutions have the same number of elements (and hence all outcome vectors have the same length), the matroid property is not used. Hence, Theorem~\ref{thm:ParetoVSlex} generalizes to all ordinal optimization problems with fixed length solutions. The following Corollary~\ref{kor:superset}, that also follows from the results in \cite{SCHAFER20}, is an immediate consequence of Theorem~\ref{thm:ParetoVSlex}.

\begin{kor}\label{kor:superset}
	The set of efficient bases of \eqref{p:mpo} is a superset of the set of efficient bases of \eqref{p:mpc} and of \eqref{p:mpcmax}. 
\end{kor}
\begin{proof}
	Theorem~\ref{thm:ParetoVSlex} implies that the efficient set of \eqref{p:mpc} can not contain any bases that are ordinally dominated w.r.t.\ $o$ since this would imply that they are also lexicographically dominated w.r.t.\ $c^{\min}$. The same argument applies to \eqref{p:mpcmax}.
\end{proof}

\begin{rem}\label{rem:reverse}
	The reverse implication of Theorem~\ref{thm:ParetoVSlex} does not hold in general, neither for $c^{\min}$ nor for $c^{\max}$. As a counter example consider the bases $B_4$ and $B_6$ from Figure~\ref{fig:allBases}. We have that $c^{\min}(B_6)=c^{\max}(B_6)=(2,1,2)$, $c^{\min}(B_4)=c^{\max}(B_4)=(1,3,1)$, $o(B_6)=(1,1,2,3,3)$ and $o(B_4)=(1,2,2,2,3)$. Hence, $c^{\min}(B_4)$ lexicographically dominates $c^{\min}(B_6)$ and $c^{\max}(B_6)$ lexicographically dominates $c^{\max}(B_4)$, while $o(B_4)$ and $o(B_6)$ are ordinally incomparable.
\end{rem}

We show in the following that in the case of matroids Corollary~\ref{kor:superset} can be strengthened. Indeed, the following result shows that the respective ordinal and lexicographic non-dominated sets are always equal and have cardinality one. This can also be observed in Example~\ref{ex:allBases} below, where all three problems \eqref{p:mpo}, \eqref{p:mpc}  and \eqref{p:mpcmax} have the same efficient and non-dominated sets.

The result can be briefly summarized as follows: Corollary~\ref{thm:ParetoVSlex} states that the efficient set of \eqref{p:mpo} is a superset of that of \eqref{p:mpc} and \eqref{p:mpcmax}. If there were two non-dominated bases $B_1$, $B_2$ for \eqref{p:mpo} and only one of them, say, basis $B_1$, was optimal for problem \eqref{p:mpc}, then the basis exchange property would imply that basis $B_2$ could be improved w.r.t.\ $c^{\min}$ by an appropriate swap operation. However, this would lead to a basis that also ordinally dominates $B_2$, contradicting the ordinal efficiency of $B_2$. This leads to the following result:

\begin{theorem}\label{thm:11}
	Let $\M_1=(\E,\I_1)$ be a matroid, let $\mathcal{X}_1\neq\emptyset$ be the set of bases of $\M_1$, and let the functions $o$, $c^{\min}$ and $c^{\max}$ be given as defined above.
	Then problems \eqref{p:mpo}, \eqref{p:mpc} and \eqref{p:mpcmax} have the same efficient set, and the corresponding non-dominated sets have cardinality one.
\end{theorem}
\begin{proof}
	We show the equality of the efficient sets of \eqref{p:mpo} and \eqref{p:mpc}. The equality of the efficient sets of \eqref{p:mpo} and \eqref{p:mpcmax} follows analogously.
	
	First observe that the non-dominated set of problem \eqref{p:mpc} has cardinality one since the lexicographical order is a total order.
	Moreover, Theorem~\ref{thm:ParetoVSlex} implies that every efficient solution of \eqref{p:mpc} is also efficient for \eqref{p:mpo}. Consequently, it is sufficient to show that all efficient solutions of \eqref{p:mpo} map to a unique non-dominated outcome vector. 
	We prove this result by contradiction. 
	
	Suppose, to the contrary, that there are two efficient bases $B_1$ and $B_2$ for \eqref{p:mpo} with $o(B_1)\neq o(B_2)$ and hence also $c^{\min}(B_1)\neq c^{\min}(B_2)$. 
	W.l.o.g.\ assume that $c^{\min}(B_1)$ lexicographically dominates $c^{\min}(B_2)$. 
	
	Let $e\in B_2\setminus B_1$ be chosen such that $o(\hat{e})\preceq o(e)$ for all $\hat{e}\in B_2\setminus B_1$, i.e., $e$ is an element of highest category among all elements in $B_2\setminus B_1$. Then the basis exchange property \eqref{basisexprop} implies that there exists an element $e'\in B_1\setminus B_2$ such that $B^*\coloneqq(B_2\cup\{e'\})\setminus\{e\} \in \X_1$, and the choice of $e$ and the fact that $c^{\min}(B_1)<_{\lex} c^{\min}(B_2)$ imply that $o(e')\preceq o(e)$. Now, if $o(e')\prec o(e)$, then $B^*$ dominates $B_2$ w.r.t.\ $o$, contradicting the assumption. Otherwise, i.e., if $o(e')= o(e)$, then $B^*$ has one more element in common with $B_1$ than $B_2$, and iterating this procedure at most $r$ times eventually yields a swap where $o(e')\prec o(e)$.
\end{proof}

\begin{kor}\label{kor:greedy}
	The ordinally non-dominated set of \eqref{p:mpo} can be computed by a greedy algorithm.
\end{kor}
\begin{proof}
	This follows immediately from Theorem~\ref{thm:11} and the matroid properties, see also \cite{hamacher1994spanning}.
\end{proof}

Note that, while Theorem~\ref{thm:11} states that the non-dominated sets of problems \eqref{p:mpo}, \eqref{p:mpc} and \eqref{p:mpcmax} have cardinality one, this does in general not transfer to the respective efficient sets. Indeed, the size of the efficient sets may grow exponentially with the problem size. As an example, consider instances with exponentially growing feasible sets and assume that all elements of $E$ are in the same ordinal category. Then, all feasible solutions of a considered problem are both ordinally and lexicographically efficient.

\section{Multi-objective Matroid Optimization}\label{sec:problem}
We extend the settings of the previous section and consider multi-objective matroid optimization problems (on a matroid $\M_1=(\E,\I_1)$ with rank $r$ and set of bases \(\X_1\)) where we combine an ordinal objective with a sum objective function with non-negative integer coefficients $w:\E\to\Zgeq$. The cost of a basis \(B\in\X_1\) w.r.t.\ this sum objective is given by $w(B)\coloneqq\sum_{e\in B} w(e)$. 

\subsection{Multi-objective Ordinal and Lexicographic Optimality and their Interrelation}\label{subsec:multiobjmodels}
If we add a sum objective function to the problems described in Section~\ref{subsec:singleobjmodels} above, we obtain the following four variants of bi- or multi-objective optimization problems involving additive as well as ordinal objective coefficients:
	
The \emph{bi-objective matroid problem with a binary objective function} \eqref{p:bbmp}
\begin{equation}\label{p:bbmp}\tag{BBMP}
	\begin{array}{ll}
		\min &\left(w(B),b(B)\right)\\
		\text{s.\,t.}&B\in \X_1 ,
	\end{array}
\end{equation}
the \emph{multi-objective matroid problem with an ordinal objective function} \eqref{p:mmpo}
\begin{equation}\label{p:mmpo}\tag{MMPO}
	\begin{array}{ll}
		\min &\left(w(B),o(B)\right)\\
		\text{s.\,t.}&B\in \X_1 ,
	\end{array}
\end{equation}
and two \emph{multi-objective matroid optimization problems with a counting objective function} \eqref{p:mmpcmax} and \eqref{p:mmpc}
\begin{equation}\label{p:mmpcmax}\tag{MMPCmax}
	\begin{array}{ll}
		\min & w(B)\\
		\lexmax &c^{\max}(B)\\
		\text{s.\,t.}&B\in \X_1
	\end{array}
\end{equation}

\begin{equation}\label{p:mmpc}\tag{MMPCmin}
	\begin{array}{ll}
		\min & w(B)\\
		\lexmin &c^{\min}(B)\\
		\text{s.\,t.}&B\in \X_1.
	\end{array}
\end{equation}

The problem~\eqref{p:bbmp} is investigated in detail in \cite{gorski10multiple} and \cite{gorski21}, where  an exhaustive swap algorithm is presented that determines a minimal complete representation of the non-dominated set (i.\,e., all non-dominated points and one efficient solution for each of them) of \eqref{p:bbmp} in polynomial time. This assumes that an oracle can determine in polynomial time if a given subset \(I\subseteq E\) is independent or not. This is, e.g., the case for graphic matroids, uniform matroids and partition matroids, see \cite{gabo:effi:84}. 

\begin{figure}
	\begin{center}
		\resizebox{\linewidth}{!}{
			\begin{tabular}{c|c|c}
				basis $B_1$ & basis $B_2$ & basis $B_3$ \\


\begin{tikzpicture}[scale=0.75]
\node[draw,circle,fill=black,scale=0.7] (1) at (0,2)[] {};
\node[draw,circle,fill=black,scale=0.7] (2) at (0,0) {};
\node[draw,circle,fill=black,scale=0.7] (3) at (1,1) {};
\node[draw,circle,fill=black,scale=0.7] (4) at (3,1) {};
\node[draw,circle,fill=black,scale=0.7] (5) at (4,2) {};
\node[draw,circle,fill=black,scale=0.7] (6) at (4,0) {};


\graph {
	(1) --["$1$",thick,draw=black!40!red] (3);
	(2) --["$2$",dashed,thick,swap, draw=yellow!40!red] (3);
	(3) --["$0$",dashed,thick,swap,draw=yellow!40!red] (4);
	(4) --["$5$",dashed,draw=yellow!40!red,thick] (5);
	(4) --["$3$",thick,swap,draw=black!40!red] (6);

};
\end{tikzpicture}

&


\begin{tikzpicture}[scale=0.75]
\node[draw,circle,fill=black,scale=0.7] (1) at (0,2)[] {};
\node[draw,circle,fill=black,scale=0.7] (2) at (0,0) {};
\node[draw,circle,fill=black,scale=0.7] (3) at (1,1) {};
\node[draw,circle,fill=black,scale=0.7] (4) at (3,1) {};
\node[draw,circle,fill=black,scale=0.7] (5) at (4,2) {};
\node[draw,circle,fill=black,scale=0.7] (6) at (4,0) {};


\graph {
(1) --["$1$",thick,draw=black!40!red] (3);
(2) --["$2$",dashed,thick,swap, draw=yellow!40!red] (3);
(3) --["$0$",dashed,thick,swap,draw=yellow!40!red] (4);
(4) --["$3$",thick,swap,draw=black!40!red] (6);
(5) --["$6$",dotted,thick,draw=black!50!green] (6);

};
\end{tikzpicture}

&


\begin{tikzpicture}[scale=0.75]
\node[draw,circle,fill=black,scale=0.7] (1) at (0,2)[] {};
\node[draw,circle,fill=black,scale=0.7] (2) at (0,0) {};
\node[draw,circle,fill=black,scale=0.7] (3) at (1,1) {};
\node[draw,circle,fill=black,scale=0.7] (4) at (3,1) {};
\node[draw,circle,fill=black,scale=0.7] (5) at (4,2) {};
\node[draw,circle,fill=black,scale=0.7] (6) at (4,0) {};


\graph {
(1) --["$1$",thick,draw=black!40!red] (3);
(1) --["$4$",dotted,thick,swap,draw=black!50!green] (2);
(3) --["$0$",dashed,thick,swap,draw=yellow!40!red] (4);
(4) --["$5$",dashed,draw=yellow!40!red,thick] (5);
(4) --["$3$",thick,swap,draw=black!40!red] (6);

};
\end{tikzpicture}

\\
				&&\\
				$w(B_1)=11$\hspace*{0.4cm} $\boldsymbol{c^{\max}(B_1)=\begin{pmatrix}
						0\\3\\2
				\end{pmatrix}}$&$w(B_2)=12$ \hspace*{0.4cm} $\boldsymbol{c^{\max}(B_2)=\begin{pmatrix}
						1\\2\\2
				\end{pmatrix}}$& $w(B_3)=13$ \hspace*{0.4cm} $c^{\max}(B_3)=\begin{pmatrix}
					1\\2\\2
				\end{pmatrix}$\\ 
				&&\\
				$\boldsymbol{o(B_1)=\begin{pmatrix}
						\oPNr\\ \oPNr \\ \oPNr \\ \rPNr \\ \rPNr
				\end{pmatrix}}$ \hspace*{0.1cm} $\boldsymbol{c^{\min}(B_1)=\begin{pmatrix}
						2\\3\\0
				\end{pmatrix}}$& $\boldsymbol{o(B_2)=\begin{pmatrix}
						\gPNr\\ \oPNr\\ \oPNr\\ \rPNr\\ \rPNr
				\end{pmatrix}}$ \hspace*{0.1cm} $\boldsymbol{c^{\min}(B_2)=\begin{pmatrix}
						2\\2\\1
				\end{pmatrix}}$ &$o(B_3)=\begin{pmatrix}
					\gPNr\\ \oPNr\\ \oPNr\\ \rPNr\\ \rPNr
				\end{pmatrix}$ \hspace*{0.1cm} $c^{\min}(B_3)=\begin{pmatrix}
					2\\2\\1
				\end{pmatrix}$ \raisebox{-1.5cm}{\rule{0pt}{1.5cm}}\\ \hline			
				basis $B_4$ & basis $B_5$ & basis $B_6$ \rule{0pt}{0.75cm}\\


\begin{tikzpicture}[scale=0.75]
\node[draw,circle,fill=black,scale=0.7] (1) at (0,2)[] {};
\node[draw,circle,fill=black,scale=0.7] (2) at (0,0) {};
\node[draw,circle,fill=black,scale=0.7] (3) at (1,1) {};
\node[draw,circle,fill=black,scale=0.7] (4) at (3,1) {};
\node[draw,circle,fill=black,scale=0.7] (5) at (4,2) {};
\node[draw,circle,fill=black,scale=0.7] (6) at (4,0) {};


\graph {
(1) --["$1$",thick,draw=black!40!red] (3);
(2) --["$2$",dashed,thick,swap, draw=yellow!40!red] (3);
(3) --["$0$",dashed,thick,swap,draw=yellow!40!red] (4);
(4) --["$5$",dashed,draw=yellow!40!red,thick] (5);
(5) --["$6$",dotted,thick,draw=black!50!green] (6);

};
\end{tikzpicture}

&


\begin{tikzpicture}[scale=0.75]
\node[draw,circle,fill=black,scale=0.7] (1) at (0,2)[] {};
\node[draw,circle,fill=black,scale=0.7] (2) at (0,0) {};
\node[draw,circle,fill=black,scale=0.7] (3) at (1,1) {};
\node[draw,circle,fill=black,scale=0.7] (4) at (3,1) {};
\node[draw,circle,fill=black,scale=0.7] (5) at (4,2) {};
\node[draw,circle,fill=black,scale=0.7] (6) at (4,0) {};


\graph {
(2) --["$2$",dashed,thick,swap, draw=yellow!40!red] (3);
(1) --["$4$",dotted,thick,swap,draw=black!50!green] (2);
(3) --["$0$",dashed,thick,swap,draw=yellow!40!red] (4);
(4) --["$5$",dashed,draw=yellow!40!red,thick] (5);
(4) --["$3$",thick,swap,draw=black!40!red] (6);

};
\end{tikzpicture}

&


\begin{tikzpicture}[scale=0.75]
\node[draw,circle,fill=black,scale=0.7] (1) at (0,2)[] {};
\node[draw,circle,fill=black,scale=0.7] (2) at (0,0) {};
\node[draw,circle,fill=black,scale=0.7] (3) at (1,1) {};
\node[draw,circle,fill=black,scale=0.7] (4) at (3,1) {};
\node[draw,circle,fill=black,scale=0.7] (5) at (4,2) {};
\node[draw,circle,fill=black,scale=0.7] (6) at (4,0) {};


\graph {
(1) --["$1$",thick,draw=black!40!red] (3);
(1) --["$4$",dotted,thick,swap,draw=black!50!green] (2);
(3) --["$0$",dashed,thick,swap,draw=yellow!40!red] (4);
(4) --["$3$",thick,swap,draw=black!40!red] (6);
(5) --["$6$",dotted,thick,draw=black!50!green] (6);

};
\end{tikzpicture}

\\
				&&\\
				$w(B_4)=14$ \hspace*{0.4cm} $c^{\max}(B_4)=\begin{pmatrix}
					1\\3\\1
				\end{pmatrix}$& $w(B_5)=14$ \hspace*{0.4cm} $c^{\max}(B_5)=\begin{pmatrix}
					1\\3\\1
				\end{pmatrix}$& $w(B_6)=14$ \hspace*{0.4cm} $\boldsymbol{c^{\max}(B_6)=\begin{pmatrix}
						2\\1\\2
				\end{pmatrix}}$\\ 
				&&\\
				$\boldsymbol{o(B_4)=\begin{pmatrix}
						\gPNr\\ \oPNr\\ \oPNr\\ \oPNr\\ \rPNr
				\end{pmatrix}}$ \hspace*{0.1cm} $\boldsymbol{c^{\min}(B_4)=\begin{pmatrix}
						1\\3\\1
				\end{pmatrix}}$&$\boldsymbol{o(B_5)=\begin{pmatrix}
						\gPNr\\ \oPNr\\ \oPNr\\ \oPNr\\ \rPNr
				\end{pmatrix}}$ \hspace*{0.1cm} $\boldsymbol{c^{\min}(B_5)=\begin{pmatrix}
						1\\3\\1
				\end{pmatrix}}$&$\boldsymbol{o(B_6)=\begin{pmatrix}
						\gPNr\\ \gPNr\\ \oPNr\\ \rPNr\\ \rPNr
				\end{pmatrix}}$ \hspace*{0.1cm} $c^{\min}(B_6)=\begin{pmatrix}
					2\\1\\2
				\end{pmatrix}$\raisebox{-1.5cm}{\rule{0pt}{1.5cm}}\\ \hline
				basis $B_7$ & basis $B_8$ & basis $B_9$ \rule{0pt}{0.75cm}\\


\begin{tikzpicture}[scale=0.75]
\node[draw,circle,fill=black,scale=0.7] (1) at (0,2)[] {};
\node[draw,circle,fill=black,scale=0.7] (2) at (0,0) {};
\node[draw,circle,fill=black,scale=0.7] (3) at (1,1) {};
\node[draw,circle,fill=black,scale=0.7] (4) at (3,1) {};
\node[draw,circle,fill=black,scale=0.7] (5) at (4,2) {};
\node[draw,circle,fill=black,scale=0.7] (6) at (4,0) {};


\graph {
(2) --["$2$",dashed,thick,swap, draw=yellow!40!red] (3);
(1) --["$4$",dotted,thick,swap,draw=black!50!green] (2);
(3) --["$0$",dashed,thick,swap,draw=yellow!40!red] (4);
(4) --["$3$",thick,swap,draw=black!40!red] (6);
(5) --["$6$",dotted,thick,draw=black!50!green] (6);

};
\end{tikzpicture}

&


\begin{tikzpicture}[scale=0.75]
\node[draw,circle,fill=black,scale=0.7] (1) at (0,2)[] {};
\node[draw,circle,fill=black,scale=0.7] (2) at (0,0) {};
\node[draw,circle,fill=black,scale=0.7] (3) at (1,1) {};
\node[draw,circle,fill=black,scale=0.7] (4) at (3,1) {};
\node[draw,circle,fill=black,scale=0.7] (5) at (4,2) {};
\node[draw,circle,fill=black,scale=0.7] (6) at (4,0) {};


\graph {
(1) --["$1$",thick,draw=black!40!red] (3);
(1) --["$4$",dotted,thick,swap,draw=black!50!green] (2);
(3) --["$0$",dashed,thick,swap,draw=yellow!40!red] (4);
(4) --["$5$",dashed,draw=yellow!40!red,thick] (5);
(5) --["$6$",dotted,thick,draw=black!50!green] (6);

};
\end{tikzpicture}

&


\begin{tikzpicture}[scale=0.75]
\node[draw,circle,fill=black,scale=0.7] (1) at (0,2)[] {};
\node[draw,circle,fill=black,scale=0.7] (2) at (0,0) {};
\node[draw,circle,fill=black,scale=0.7] (3) at (1,1) {};
\node[draw,circle,fill=black,scale=0.7] (4) at (3,1) {};
\node[draw,circle,fill=black,scale=0.7] (5) at (4,2) {};
\node[draw,circle,fill=black,scale=0.7] (6) at (4,0) {};


\graph {
(2) --["$2$",dashed,thick,swap, draw=yellow!40!red] (3);
(1) --["$4$",dotted,thick,swap,draw=black!50!green] (2);
(3) --["$0$",dashed,thick,swap,draw=yellow!40!red] (4);
(4) --["$5$",dashed,draw=yellow!40!red,thick] (5);
(5) --["$6$",dotted,thick,draw=black!50!green] (6);

};
\end{tikzpicture}

\\
				&&\\
				$w(B_7)=15$ \hspace*{0.4cm} $\boldsymbol{c^{\max}(B_7)=\begin{pmatrix}
						2\\2\\1
				\end{pmatrix}}$& $w(B_8)=16$ \hspace*{0.4cm}  $c^{\max}(B_8)=\begin{pmatrix}
					2\\2\\1
				\end{pmatrix}$& $w(B_9)=17$ \hspace*{0.4cm}  $\boldsymbol{c^{\max}(B_9)=\begin{pmatrix}
						2\\3\\0
				\end{pmatrix}}$\\ 
				&&\\
				$\boldsymbol{o(B_7)=\begin{pmatrix}
						\gPNr\\ \gPNr\\ \oPNr\\ \oPNr\\ \rPNr
				\end{pmatrix}}$ \hspace*{0.1cm} $\boldsymbol{c^{\min}(B_7)=\begin{pmatrix}
						1\\2\\2
				\end{pmatrix}}$ &$o(B_8)=\begin{pmatrix}
					\gPNr\\ \gPNr\\ \oPNr\\ \oPNr\\ \rPNr
				\end{pmatrix}$ \hspace*{0.1cm} $c^{\min}(B_8)=\begin{pmatrix}
					1\\2\\2
				\end{pmatrix}$&$\boldsymbol{o(B_9)=\begin{pmatrix}
						\gPNr\\ \gPNr\\ \oPNr\\ \oPNr\\ \oPNr
				\end{pmatrix}}$ \hspace*{0.1cm} $\boldsymbol{c^{\min}(B_9)=\begin{pmatrix}
						0\\3\\2
				\end{pmatrix}}$
		\end{tabular}	}
	\end{center}
	\caption{All bases of the graphic matroid introduced in Figure~\ref{fig:bbmpVsMop:2} together with the objective values $w$, $o$, $c^{\min}$ and $c^{\max}$, where we write $\gPNr$ for green-dotted, $\oPNr$ for orange-dashed, and $\rPNr$ for red-solid edges. When only considering the sum objective $w$, then $B_1$ is optimal, and when only considering the objective functions $o$, $c^{\min}$ or $c^{\max}$, respectively, then $B_9$ is the unique efficient basis. For the problems \eqref{p:mmpo}, \eqref{p:mmpc} and \eqref{p:mmpcmax} that combine $w$ with $o$, $c^{\min}$ and $c^{\max}$, respectively, the non-dominated outcome vectors are indicated by printing the partial objective vectors $o$, $c^{\min}$ and $c^{\max}$ in bold.}
	\label{fig:allBases}
\end{figure}
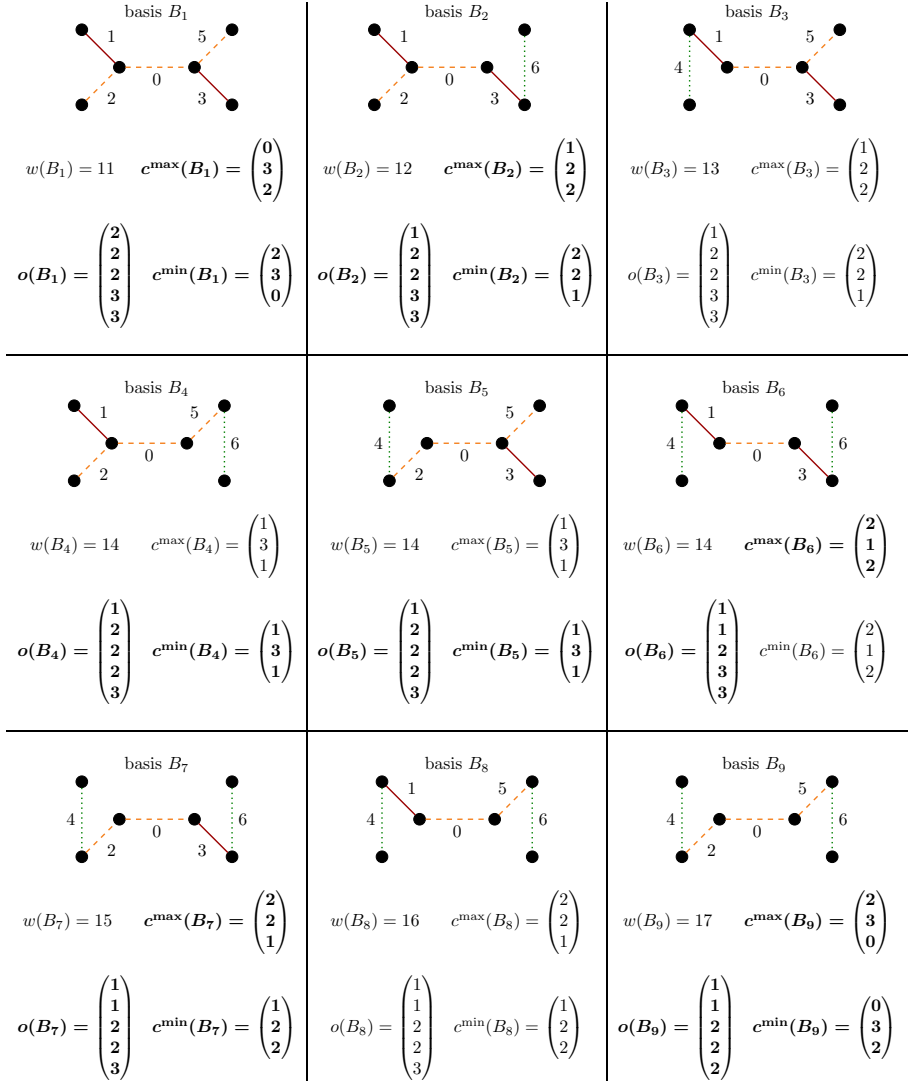

The similarities and differences between the problems \eqref{p:mmpo}, \eqref{p:mmpcmax} and \eqref{p:mmpc} are illustrated at the following example of a graphic matroid:

\begin{example}\label{ex:allBases}
	Consider the graphic matroid introduced in Figure~\ref{fig:bbmpVsMop:2}. 
	Its bases are enumerated and illustrated with their weight functions $w$, $o$, $c^{\min}$ and $c^{\max}$ in Figure~\ref{fig:allBases}. It is easy to see that, in accordance with Theorem~\ref{thm:11}, the unique efficient solution w.r.t.\ all of the individual objective functions $o$, $c^{\min}$ and $c^{\max}$ is the basis $B_{9}$. 
	
	The corresponding multi-objective problems that additionally consider the sum objective function $w$ all have larger non-dominated sets in this example. The respective non-dominated outcome vectors of the multi-objective problems that combine $w$ with the objective functions $o$, $c^{\min}$ and $c^{\max}$, respectively, are highlighted in Figure~\ref{fig:allBases} by printing the latter components, i.e., $o$, $c^{\min}$ and $c^{\max}$, in bold. Note that the basis $B_1$ is efficient in all three cases since it is the unique minimizer of $w$.
\end{example}

\subsection{Interrelation Between \eqref{p:mmpo}, \eqref{p:mmpc} and \eqref{p:mmpcmax}}
When moving from optimization problems with only one ordinal objective function to multi-objective problems that additionally include a sum objective $w$, as in the multi-objective problems \eqref{p:mmpo}, \eqref{p:mmpc} and \eqref{p:mmpcmax}, the situation is much more complex than that described in Section~\ref{sec:inerrelateo} above. Indeed, while Corollary~\ref{kor:superset} can be adapted to the new situation, Theorem~\ref{thm:11} does not transfer to the multi-objective case. A corresponding counter example will be given below.

\begin{theorem}\label{cor:superset2}
		The set of efficient bases of \eqref{p:mmpo} is a superset of the set of efficient bases of \eqref{p:mmpc} and of \eqref{p:mmpcmax}.
\end{theorem}
	\begin{proof}
		We prove the result for $c^{\min}$. The corresponding result for $c^{\max}$ follows analogously, noting that \eqref{p:mmpc} involves lexicographic minimization while \eqref{p:mmpcmax} involves lexicographic maximization.
		
		We prove this result by contradiction. Hence, let $\bar{B}$ be an efficient basis for \eqref{p:mmpc} but not for \eqref{p:mmpo}. Then there exists a basis $B^*$ with $w(B^*)\leq w(\bar{B})$, $o(B^*) \preceqq o(\bar{B})$, and $(w(B^*),o(B^*)) \neq (w(\bar{B}), o(\bar{B}))$. First note that $o(B^*) \preceqq o(\bar{B})$ implies that $c^{\min}(B^*)\leqq_{\lex} c^{\min}(\bar{B})$, by Theorem~\ref{thm:ParetoVSlex}. We distinguish two cases: Either $w(B^*)< w(\bar{B})$ and $c^{\min}(B^*)\leqq_{\lex} c^{\min}(\bar{B})$, or $w(B^*)= w(\bar{B})$ and $c^{\min}(B^*)<_{\lex} c^{\min}(\bar{B})$. However, both cases are in contradiction with the efficiency of $\bar{B}$ for problem \eqref{p:mmpc}.
		\end{proof}
	
However, as was to be expected, Theorem~\ref{thm:11} does not generalize to the multi-objective case as is shown by the following counter example:

\begin{example}\label{ex:mocase}
	Consider again the graphic matroid introduced in Example~\ref{ex:allBases} and the set of all of its bases illustrated in Figure~\ref{fig:allBases}.
	The efficient bases for  problem \eqref{p:mmpo} are the bases $B_1,B_2,B_4,B_5,B_6,B_7,B_9$,  
	while the efficient bases for the problem \eqref{p:mmpc} are given by $B_1,B_2,B_4,B_5,B_7,B_9$, 
		and the efficient bases for problem \eqref{p:mmpcmax} are given by
		$B_1,B_2,B_6,B_7,B_9$.
	Hence, basis $B_6$ is efficient for \eqref{p:mmpo} but not for \eqref{p:mmpc}, 
	and the two bases $B_4$ and $B_5$ are efficient for \eqref{p:mmpc} but not for \eqref{p:mmpcmax}.
	Thus, Theorem~\ref{thm:11} does not generalize to the multi-objective problems \eqref{p:mmpo}, \eqref{p:mmpc} and \eqref{p:mmpcmax}.
\end{example}

One could conjecture from Example~\ref{ex:allBases} that every efficient basis for \eqref{p:mmpo} is efficient for at least for one of the problems \eqref{p:mmpc} or \eqref{p:mmpcmax}. However, this also does not hold in general as the following example shows.

\begin{figure}[htb]
	\begin{center}


\begin{tikzpicture}[scale=0.8]
\node[draw,circle,thick,fill=black!10] (1) at (0,0)[] {};
\node[draw,circle,thick,fill=black!10] (2) at (2,0) {};
\node[draw,circle,thick,fill=black!10] (3) at (4,0) {};
\node[draw,circle,thick,fill=black!10] (4) at (6,0) {};
\node[draw,circle,thick,fill=black!10] (5) at (8,0) {};


\graph {
(1) --["\footnotesize $0$",thick,draw=black!40!red, bend right=60] (2);
(1) --["\footnotesize $1$",thick,draw=yellow!40!red, dashed] (2);
(1) --["\footnotesize $2$",thick,draw=black!50!green, dotted, bend left=60] (2);
         
(2) --["\footnotesize $0$",thick,draw=black!40!red, bend right=60] (3);
(2) --["\footnotesize $1$",thick,draw=yellow!40!red, dashed] (3);
(2) --["\footnotesize $2$",thick,draw=black!50!green, dotted, bend left=60] (3);
         
(3) --["\footnotesize $0$",thick,draw=black!40!red, bend right=60] (4);
(3) --["\footnotesize $1$",thick,draw=yellow!40!red, dashed] (4);
(3) --["\footnotesize $2$",thick,draw=black!50!green, dotted, bend left=60] (4);
         
(4) --["\footnotesize $0$",thick,draw=black!40!red, bend right=60] (5);
(4) --["\footnotesize $1$",thick,draw=yellow!40!red, dashed] (5);
(4) --["\footnotesize $2$",thick,draw=black!50!green, dotted, bend left=60] (5);

};
\end{tikzpicture}

		\vspace*{.5cm}
		\begin{tabular}{c|c|cc|cc|ccc}
			$w=0$&$w=1$&\multicolumn{2}{c|}{$w=2$}&\multicolumn{2}{c|}{$w=3$}&\multicolumn{3}{c}{$w=4$}\\
			\hline 
			$\begin{pmatrix}3\\ 3\\ 3\\ 3\\ \end{pmatrix}$\phantom{\rule{0pt}{3.1em}} &
			$\begin{pmatrix}2\\ 3\\ 3\\ 3\\ \end{pmatrix}$ &
			$\begin{pmatrix}1\\ 3\\ 3\\ 3\\ \end{pmatrix}$ &
			$\begin{pmatrix}2\\ 2\\ 3\\ 3\\ \end{pmatrix}$ &
			$\begin{pmatrix}1\\ 2\\ 3\\ 3\\ \end{pmatrix}$ &
			$\begin{pmatrix}2\\ 2\\ 2\\ 3\\ \end{pmatrix}$ &
			$\begin{pmatrix}1\\ 1\\ 3\\ 3\\ \end{pmatrix}$ &
			$\begin{pmatrix}1\\ 2\\ 2\\ 3\\ \end{pmatrix}$ &
			$\begin{pmatrix}2\\ 2\\ 2\\ 2\\ \end{pmatrix}$ 
		\end{tabular}
	\end{center}
	\caption{All possible outcome vectors $o(B)$ with $w(B)\in\{0,\dots,4\}$ for a graphic matroid with non-negative integer-valued costs $w$ and three categories (1:green-dotted, 2:orange-dashed and 3:red-solid). }
	\label{fig:K3}
\end{figure}
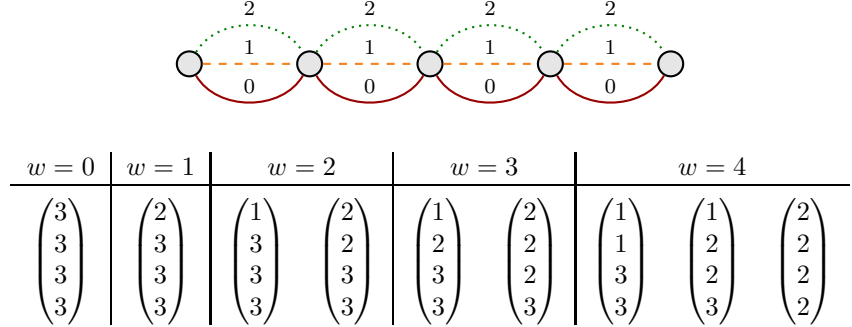

\begin{example}
	Consider the graphic matroid shown in Figure~\ref{fig:K3}. We focus on all bases $B\in\X_1$ that have an objective value of $w(B)=4$ in the sum objective. Note that these bases can only be dominated by other bases $\hat{B}$ with $w(\hat{B})\leq w(B)$, and hence we restrict our analysis on those bases in Figure~\ref{fig:K3}. First observe that all bases $B\in\X_1$ with $w(B)=4$ map to one of the three possible outcome vectors $o(B)\in\{(1,1,3,3), (1,2,2,3), (2,2,2,2)\}$, which are all non-dominated for \eqref{p:mmpo}. Their corresponding counting vectors $c^{\min}$ are $(2,0,2)$, $(1,2,1)$ and $(0,4,0)$, where the last one is the only one that is lexicographically non-dominated. For $c^{max}$ the counting vectors are the same, but the first one is lexicographically non-dominated. Consequently, the counting vector $(2,0,2)$ is neither lexicographically non-dominated for \eqref{p:mmpc} nor for \eqref{p:mmpcmax}, but it is non-dominated for \eqref{p:mmpo}. 
\end{example}

\section{Matroid Intersection for Ordinal Constraints}
\label{sec:matroidInter}
In the following we show that the three problems \eqref{p:mmpo}, \eqref{p:mmpc} and \eqref{p:mmpcmax} can be solved using a series of matroid intersection problems. The approach is based on variants of $\varepsilon$-constraint scalarizations of problem \eqref{p:mmpo} with appropriately selected optimization objective and constraints. Furthermore, we show that matroid intersection problems can be used to solve even problems with several ordinal objective functions and one sum objective.

\subsection{Variants of $\eps$-Constraint Scalarizations}\label{subsec:eps}
We consider an equality-constrained scalarization of \eqref{p:mmpo} (where equality constraints are used rather than inequality constraints as is commonly the case in $\eps$-constraint scalarizations), given by
\begin{equation}\label{eq:epseq}
	\begin{array}{ll}
		\min&w(B)\\
		\text{s.\,t.}& o_i(B) = \eps_i, \quad i=1,\dots,r\\
		& B \in \X_1
	\end{array}
\end{equation}
with right-hand side vector $\eps\in \mathcal{C}^r$. 
Intuitively, problem \eqref{eq:epseq} specifies \emph{exactly} how many elements of each category must be chosen, and hence each feasible basis $B\in\X_1$ of \eqref{eq:epseq} maps to the same ordinal vector $o(B)$. 
Depending on the choice of $\eps$, problem \eqref{eq:epseq} may be infeasible (if there is no $B\in\X_1$ with $o(B)=\eps$), yield an efficient solution $B^*$ for \eqref{p:mmpo} (if there is \emph{no} $B\in\X_1$ with $w(B)=w(B^*)$ and $o(B)\preccurlyeq o(B^*)$), or yield a dominated solution $\hat{B}$ for \eqref{p:mmpo}  (if there \emph{is} a $B\in\X_1$ with $w(B)=w(\hat{B})$ and $o(B)\preccurlyeq o(\hat{B})$). Note that suitable choices for $\eps$ satisfy $\eps_1\preceq\cdots\preceq\eps_r$ since the components of $o(B)$ are always in non-decreasing order and hence problem \eqref{eq:epseq} is certainly infeasible otherwise. In the following, we denote all such \emph{suitable right-hand-side vectors} by  $\Upsilon\coloneqq\{\eps\in \mathcal{C}^r \colon \eps_1\preceq\cdots\preceq\eps_r\}$.
 
Since problem \eqref{eq:epseq} can be interpreted as a variant of the ``classical'' $\eps$-constraint scalarization in multi-objective optimization, see, e.g., \cite{Ehrg05}, the following result is not surprising and follows basically by the same arguments.

\begin{theorem}\label{th:econsteqforo}
	The non-dominated set of problem \eqref{p:mmpo} can be determined by solving problem \eqref{eq:epseq} for all suitable right-hand side vectors  $\eps\in \Upsilon$ and filtering out all dominated outcome vectors.
\end{theorem}

\begin{proof}
	Let $(w(B^*),o(B^*))$ be a non-dominated outcome vector for \eqref{p:mmpo} with pre-image $B^*\in\X_1$. Then $B^*$ is optimal for problem \eqref{eq:epseq} with $\eps\coloneqq o(B^*)\in\Upsilon$. Thus, every non-dominated outcome vector of problem \eqref{p:mmpo} can be determined by solving an appropriate scalarization \eqref{eq:epseq}. The non-dominated set is then obtained by employing a dominance filtering to the set of all obtained outcome vectors.
\end{proof}

Now let a suitable constraint vector $\eps\in \Upsilon$ be given, i.e., $\eps$ satisfies $\eps_1\preceq\cdots\preceq\eps_r$. Then we can define an associated suitable counting vector $u\in\Zgeq^K$ by setting $u_i\coloneqq\vert\{j\in\{1,\dots,r\}:\eps_j=\eta_i\}\vert$ for all $i=1,\dots,K$, where, by definition, we have that $\sum_{i=1}^K u_i=r$. We denote by $U\coloneqq\{u\in\Zgeq^K\colon \sum_{i=1}^K u_i=r\}$ the set of all \emph{suitable counting vectors}.

\begin{lem}\label{kor:epsversusu}
 There is a one-to-one correspondence between suitable right-hand-side vectors $\eps\in \Upsilon$ and suitable counting vectors $u\in U$.
\end{lem}
\begin{proof}
First consider the case that a suitable constraint vector $\eps\in\Upsilon$ is given. Then an associated suitable counting vector $u\in U$ can be determined from $\eps$ as described above, i.e., by setting  $u_i\coloneqq\vert\{j\in\{1,\dots,r\}:\eps_j=\eta_i\}\vert$ for all $i=1,\dots,K$. Conversely, if a suitable counting vector $u\in U$ is given, then we can determine associated suitable values for $\eps\in\Upsilon$ by setting
$\eps_j\coloneqq\eta_i$, where the ordinal level $i\in\{1,\dots,K\}$ is chosen such that $\sum_{l=1}^{i}u_l\leq j$ and $\sum_{l=1}^{i-1}u_l> j$, for all $j=1,\dots,r$.
\end{proof}

Lemma~\ref{kor:epsversusu} implies that problem \eqref{eq:epseq} can be equivalently written as
\begin{equation}\label{eq:ueq}
		\begin{array}{ll}
			\min&w(B)\\
			\text{s.\,t.}& c_i(B) = u_i, \quad i=1,\dots,K\\
			& B \in \X_1,
		\end{array}
\end{equation}
where the right-hand side vector $u\in U$ is chosen as a suitable counting vector, i.e., $u\in\Zgeq^K$ and $\sum_{i=1}^K u_i=r$, and $c$ is the counting objective introduced in Section~\ref{subsec:singleobjmodels}.
Moreover, since $\sum_{i=1}^K c_i(B)=r=\sum_{i=1}^K u_i$ for all feasible bases $B\in\X_1$, the equality constraints in \eqref{eq:ueq} can be replaced by inequality constraints without changing the feasible set. Problem~\eqref{eq:ueq} is thus equivalent to the following variant of $\eps$-constraint scalarization that relates to problem \eqref{p:mmpc} 

	\begin{equation}\label{eq:eps-cmin}
		\begin{array}{ll}
			\min&w(B)\\
			\text{s.\,t.}& c^{\min}_i(B)\leq u_{K-i+1}, \quad i=1,\dots,K\\
			& B \in \X_1.
		\end{array}
\end{equation}

\begin{kor}\label{kor:econsteqforo}
	The non-dominated set of problem \eqref{p:mmpo} can be determined by solving problem \eqref{eq:ueq} (or problem \eqref{eq:eps-cmin}) for all suitable counting vectors $u\in U$, and filtering out all dominated outcome vectors. The non-dominated sets of problems \eqref{p:mmpc} and \eqref{p:mmpcmax} can be obtained from this set by further filtering out all lexicographically dominated outcome vectors.
\end{kor}
\begin{proof}
	The result follows immediately from Theorems~\ref{cor:superset2} and~\ref{th:econsteqforo}, using the equivalence of the formulations \eqref{eq:epseq}, \eqref{eq:ueq} and \eqref{eq:eps-cmin}.
\end{proof}

We emphasize that problem \eqref{eq:eps-cmin} remains meaningful when a non-suitable counting vector $u\in\Zgeq^K$ with $\sum_{i=1}^K u_i > r$ is used as right-hand-side vector. Indeed, when considering a suitable counting vector $u\in U$ and a non-suitable counting vector $\hat{u}\geqslant u$, then $\hat{u}$ yields a relaxation of problem \eqref{eq:eps-cmin} with  $u$ as the right-hand-side vector. Nevertheless, the constraint $B\in\X_1$ guarantees that only bases of $\M_1$ are returned, and hence $\sum_{i=1}^K c_i^{\min}(B)=r$ remains satisfied also in this case. Moreover, using $\bar{u}\coloneqq(r,\dots,r)\in\Zgeq^K$ yields a ``complete'' relaxation in the sense that constraints $c^{\min}_i(B)\leq \bar{u}_{K-i+1}=r$, $i=1,\dots,K$, are satisfied for \emph{all} bases $B\in\X_1$ and hence redundant in this case.

Using the above results, the cardinality of the non-dominated set of problem \eqref{p:mmpo} (and hence also of problems \eqref{p:mmpc} and \eqref{p:mmpcmax}) can be bounded. Note that this is in analogy to the results obtained in \cite{SCHAFER20} and \cite{SCHAFER:knapsack} for ordinal shortest path and ordinal knapsack problems, respectively.
Indeed, the number of equality-constraint scalarizations \eqref{eq:ueq} that need to be solved in order to guarantee that all non-dominated outcome vectors of problem \eqref{p:mmpo} are found is polynomially bounded. This can be seen from the fact that 
the number of suitable counting vectors $u\in U$, i.e., the number of $K$-dimensional non-negative integer vectors that satisfy $\sum_{j=1}^K u_j=r$, is given by $\binom{r+K-1}{K-1}=\bigo(r^{K-1})$ (assuming that $K$ is constant), i.\,e., it is equal to the number of multisets of cardinality $K-1$ taken from a set of size $r+1$. This number is also known as occupancy number, see, e.g., \cite{Feller68}.
We obtain the following result.

\begin{theorem}\label{cor:poly}
		The cardinality of the non-dominated set of problem \eqref{p:mmpo} is bounded by $\bigo(r^{K-1})$, which is polynomial in $r$ as long as $K$ is constant.
\end{theorem}

\subsection{Matroid Intersection}
We focus on the $\eps$-constraint scalarization-variant \eqref{eq:eps-cmin} in the following and show that it can be equivalently formulated as a matroid intersection problem. Towards this end, let an arbitrary but fixed, suitable counting vector $u\in U$ be given as right-hand-side vector in problem \eqref{eq:eps-cmin}.

Now consider the partition $\E=E_1\cup E_2\cup\ldots\cup E_K$ of the ground set $E$ of $\M_1$, where $E_j\coloneqq \{e\in \E \colon o(e)=\eta_j \}$ for $j=1,\dots,K$, i.e., $E_j$ contains all elements from $E$ that are in category $\eta_j$. Given this partition of $E$, let $\M_2(u)=(\E,\J(u))$ be an associated partition matroid with independent sets given by   
$\J(u)\coloneqq \{J\subseteq \E \colon\vert J\cap E_j\vert\leq u_j, \; 1\leq j\leq K\}$.
Then problem \eqref{eq:eps-cmin} can be solved using the matroid intersection problem
\begin{equation}\label{p:mi-indep}
	\begin{array}{ll}
		\min&w(I)\\
		\text{s.\,t.}& I\in \I_1\cap \I_2(u)\\
		& \vert I\vert=\max\{\vert J\vert:J\in\I_1\cap\I_2(u)\}.
	\end{array}
\end{equation}
Note that the second constraint in \eqref{p:mi-indep} is needed since otherwise, $B=\emptyset$ would always be optimal.

\begin{theorem} Let $u\in \Zgeq^K$ be arbitrary but fixed. 
 If problem \eqref{eq:eps-cmin} is feasible, then problems \eqref{eq:eps-cmin} and \eqref{p:mi-indep} are equivalent. Moreover, if problem \eqref{eq:eps-cmin} is infeasible, then every optimal solution $B^*$ of problem \eqref{p:mi-indep} satisfies $|B^*|<r$.
\end{theorem}

\begin{proof}
 We first show that when problem \eqref{eq:eps-cmin} is feasible for the given suitable counting vector $u$, then  problems \eqref{eq:eps-cmin} and \eqref{p:mi-indep} have the same feasible sets. Indeed, in this case there exists a basis $\hat{B}\in\X_1$ that satisfies $c(\hat{B})\leqq u$, and hence $\hat{I}\coloneqq \hat{B}$ with $|\hat{I}|=r$ is feasible for \eqref{p:mi-indep}. This implies that all feasible solutions of \eqref{p:mi-indep} have cardinality $r$ and are thus bases of $\M_1$. In this situation, the constraints $c^{\min}_i(B)\leq u_{K-i+1}$, $i=1,\dots,K$ (for \eqref{eq:eps-cmin}) and $I\in\I_2(u)$ (for \eqref{p:mi-indep}) are equivalent.  Since both problems also have the same objective function, they are clearly equivalent in this case. If, however, problem~\eqref{eq:eps-cmin} is infeasible for the current choice of $u$, then the matroid intersection problem \eqref{p:mi-indep} is still feasible, but returns an optimal solution $B^*$ with $|B^*|<r$. This situation can be easily recognized.
\end{proof}

The advantage of this reformulation is that the matroid intersection problem \eqref{p:mi-indep} can be solved by the polynomial time \emph{matroid intersection algorithm} (MI) of \cite{edmonds71matroids}. We refer to \cite{Schrijver13} for a proof of its correctness and of its polynomial run time. Note that, while the matroid intersection algorithm is originally formulated for maximization problems, it can also be applied to problem \eqref{p:mi-indep} by multiplying all weights $w(e)$ by $-1$, $e\in E$, and hence maximizing $-w(e)$ rather than minimizing $w(e)$ in \eqref{p:mi-indep}.  
In our implementation we use a variant of the Floyd-Warshall algorithm (see, e.g., \cite{ahuja1993network}) to realize the required shortest path computations. Note that, alternatively, the weights could be further transformed such that non-negative weight coefficients are obtained. Then, the algorithm of Dijkstra could be applied for the relevant distance computations, see \cite{FRANK1981328} and \cite{brezovec1986two} for more details.

Note that the cardinality constraint in the general formulation \eqref{p:mi-indep} of the matroid intersection problem  can be omitted if $\I_1$ is replaced by $\X_1$, i.e., the set of all bases of the matroid $\M_1$. Hence, we can alternatively solve the problem
\begin{equation}\label{p:mi-basis}
	\begin{array}{ll}
		\min&w(B)\\
		\text{s.\,t.}& B\in \X_1\cap \I_2(u).
	\end{array}
\end{equation}
This is realized in Algorithm~\ref{alg:mioc} by only considering optimal solutions of problem~\eqref{p:mi-indep} that are actually bases of $\M_1$, see lines 4 and 5 in Algorithm~\ref{alg:mioc} below.

\subsection{Algorithmic Consequences}\label{subsec:alg}

Theorems~\ref{th:econsteqforo} and~\ref{cor:poly} imply that all solutions of \eqref{p:mmpo} can be determined by a polynomial number of matroid intersections applied on \eqref{p:mi-indep}. The structure of this procedure is given in Algorithm~\ref{alg:mioc}. Note that this algorithm can be easily adapted to solve \eqref{p:mmpc} or \eqref{p:mmpcmax}. Since the non-dominated set of \eqref{p:mmpo} is a superset of the corresponding non-dominated sets of \eqref{p:mmpc} and \eqref{p:mmpcmax}, only a slight modification of the filtering step in line 6 is necessary.

\begin{algorithm}
	\SetAlgoLined 
	\LinesNumbered
	\KwIn{Matroid $\M_1=(E,\I_1)$, sum objective function $w$ and ordinal objective function $o$}
	\KwOut{Non-dominated set of problem \eqref{p:mmpo}}
	\caption{Matroid Intersection for Ordinal Constraints ($\mioc$($\M_1,w,o$))}\label{alg:mioc}
	$X\coloneqq\emptyset$\;
	\ForEach{\(u\in U\)}{
		Solve \eqref{p:mi-indep} with (MI) and save the obtained independent set $I^*$\;
		\If{$\vert I^*\vert=r$}{Set $X=X\cup\{I^*\}$\;}
	}
	Filter the efficient independent sets of \(X\) w.r.t.\ \eqref{p:mmpo} and save the corresponding outcome
	vectors in $\outNd^{O}$ \;
	\Return $\outNd^{O}$
\end{algorithm}

It is possible to improve the performance of Algorithm~\ref{alg:mioc} by reducing the number of considered bounds \(u\in U\), i.e., the number of solved matroid intersections. This can be achieved by initially solving \eqref{p:mi-indep} with $u=(r,...,r)\in\mathbb{R}^K$, which returns a weakly efficient basis $B^*$ (assuming $\X_1\neq\emptyset$) with the smallest possible cost $w^*$. Consequently, only upper bounds \(u\in U\) such that $(u_K,\dots,u_1)\leq_{\lex} c^{\min}(B^*)$ 
have to be considered. Thus, we modify lines 1--2 in Algorithm~\ref{alg:mioc} accordingly and obtain Algorithm~\ref{alg:mioc-mmpo}.

\begin{algorithm}
	\SetAlgoLined 
	\LinesNumbered
	\KwIn{Matroid $\M_1=(E,\I_1)$, sum objective function $w$ and  ordinal objective function $o$}
	\KwOut{Non-dominated set of problem \eqref{p:mmpo}}
	\caption{Improved Initialization of Matroid Intersection for Ordinal Constraints 
	($\mioc^O$($\M_1,w,o$))}\label{alg:mioc-mmpo}
	 $u=(r,...,r)$\;
	Solve \eqref{p:mi-indep} with (MI) and save the obtained basis $B^*$ \;
	Set $X=\{B^*\}$\;
	\ForEach{$u\in \{v\in U \colon (v_K,\dots,v_1)\leq_{\lex}c^{\min}(B^*)\}$ }{
	\textit{run lines 3--5 of Algorithm~\ref{alg:mioc}}
	}
	Filter the efficient independent sets of \(X\) w.r.t.\ \eqref{p:mmpo} and save the corresponding outcome
	vectors in $\outNd^{O}$\;
	\Return $\outNd^{O}$
\end{algorithm}

Note that, in the worst case, this initialization yields no reduction of the running time, since there might exists a basis $B'$ that minimizes $w$ and for which $c^{\min}(B')=(r,0,\dots,0)$. However, in our numerical tests this procedure often leads to a significant reduction of the number of iterations, as described in Section~\ref{sec:numeric}.
Note that the initial bound in Algorithm~\ref{alg:mioc-mmpo} is not a suitable counting vector as defined in Section~\ref{subsec:eps}, i.e., \(u=(r,\ldots,r)\notin U\). Since we consider in the following often relaxations of suitable subproblems, we use the notation \(\bar{U}\coloneqq\{u\in\Zgeq^K\colon \sum_{i=1}^K u_i \geq r,\; u_i\leq r,\; i=1,\ldots,K\}\) to denote the considered upper bound set.

Based on the fact that the lexicographic order is a total order, Algorithm~\ref{alg:mioc-mmpo} can be further improved when applied on problem~\eqref{p:mmpc}. In this case, the considered upper bound set can also be reduced during the course of the algorithm. 
The initialization of the bound set \(U\) is analogous to Algorithm~\ref{alg:mioc-mmpo}, i.\,e., we solve the matroid intersection problem for $u=(r,...,r)\in\mathbb{R}^K$. Let $B^*$ be the obtained weakly efficient basis of \eqref{p:mmpc} minimizing the sum objective function $w$.
Then, it is sufficient to solve subproblems with upper bounds \(u\in \bar{U}\) such that $(u_K,\dots,u_1)\leq_{\lex} c^{\min}(B^*)$. Due to the lexicographic order we can explicitly enumerate  the new upper bounds $u$ to be considered as  $(u_K,\dots,u_1)\in\{(c_1^{\min}(B^*)-1,r,\dots,r),(c_1^{\min}(B^*),c_2^{\min}(B^*)-1,r,\dots,r),\dots,(c_1^{\min}(B^*),\dots,c_{K-2}^{\min}(B^*),c_{K-1}^{\min}(B^*)-1,r)\}$ such that $u\geqslant0$. These upper bounds are added to the list $\bar{U}$ of open subproblems and sorted in lexicographically increasing order.
Whenever a new candidate for an efficient basis is found, we update the list of open subproblems $\bar{U}$ and re-sort it.
In Algorithm~\ref{alg:mioc-mmpc} this procedure is repeated until $\bar{U}=\emptyset$. To simplify the notation we slightly abuse the notation and consider \(\bar{U}\) to be a sorted list, referring by \(\bar{U}[1]\) to the first element of this list.
Note that an analogous solution algorithm can be formulated for the corresponding lexicographic maximization problem~\eqref{p:mmpcmax}.

\begin{algorithm}
	\SetAlgoLined 
	\LinesNumbered
	\KwIn{Matroid $\M_1=(E,\I_1)$, sum objective function $w$ and an ordinal objective function $o$}
	\KwOut{Non-dominated set of problem \eqref{p:mmpc}}
	\caption{Matroid Intersection for Ordinal Constraints for \eqref{p:mmpc} ($\mioc^{C_{\min}}$($\M_1,w,o$))}\label{alg:mioc-mmpc}
	$\bar{U}\coloneqq\{u\in\Zgeq^K\colon \sum_{i=1}^K u_i \geq r\}$, $u=(r,...,r)\in\Zgeq^K$\;
	Solve \eqref{p:mi-indep} with (MI) and save the obtained basis $B^*$ \;
	Set $X= \{B^*\}$\;
	$\bar{U}\coloneqq \{u\in\bar{U}\colon (u_K,\dots,u_1)\leq_{\lex}c^{\min}(B^*)\}$, sort \(U\) in lexicographically increasing order \;
	\While{$\bar{U}\neq \emptyset$}{
		\(u\coloneqq\bar{U}[1]\)\tcp*{pop lexicographically smallest bound} 
		\(\bar{U}\coloneqq\bar{U}[2,\ldots,\text{end}]\)\;
		Solve \eqref{p:mi-indep} with (MI) and save the obtained independent set $I^*$\;
		\If{$\vert I^*\vert=r$}{Set $X=X\cup\{I^*\}$\;
		$\bar{U}\coloneqq \{u\in\bar{U}\colon (u_K,\dots,u_1)\leq_{\lex}c^{\min}(I^*)\}$, sort \(U\) in lexicographically increasing order \;}		
	}
	Filter the efficient independent sets of $X$ with respect to the problem \eqref{p:mmpc} and save the corresponding outcome vectors in $\outNd^{C_{\min}}$ \;
	\Return $\outNd^{C_{\min}}$
\end{algorithm}

\subsection{Problems with Several Ordinal Objective Functions}
Problem~\eqref{p:mmpo} can be generalized by considering \(p\geq 2\) objective functions with ordinal weights.
This can be illustrated at a graph whose edges are classified w.r.t.\  two types of categories, for example, colors (e.g., green, orange, red) and letters (e.g., A, B). Then every edge is in exactly one of the following categories: green-A, green-B, orange-A, orange-B, red-A or red-B. 
These combinations of categories are a-priori not completely ordered, since, in general, neither green-B is preferred over red-A, nor red-A is preferred over green-B. However, such problems can be considered in the context of combined orderings as introduced in Section~\ref{subsec:multiobjective}.

Multi-objective matroid problems with one sum objective function and several ordinal objective functions can be handled analogously to multi-objective matroid problems with one ordinal objective function \eqref{p:mmpo}. Without going much into detail, we shortly describe the formulation of an associated weighted matroid intersection problem that generalizes problem \eqref{p:mi-indep}.

Let $p$ denote the number of ordinal objective functions $o^i$, $i=1,\dots,p$, let $K_{i}$ denote the number of categories for the $i$-th ordinal objective, and let $\eta_{ij}\in\mathcal{C}$ denote the $j$-th category of the $i$-th ordinal objective, $j=1,\ldots,K_i$, where $\eta_{ij}\prec\eta_{ik}$ whenever $j<k$. Then we can define a partition matroid $\M_3$ (generalizing $\M_2$) by partitioning the ground set  $\E=\bigcup_{i=1}^p\bigcup_{j=1}^{K_{i}} E_{ij}$, where $E_{ij}\coloneqq \{e\in \E:o^i(e)=\eta_{ij}\}$ for $i=1,\dots,p$ and $j=1,\dots,{K_{i}}$. The set of independent sets of $\M_3$ is given by $\I_3\coloneqq \{J\subseteq \E:\vert J\cap E_{ij}\vert\leq u_{ij},\; 1\leq i\leq p\text{ and } j\in\{1,\ldots, {K_{i}}\}\}$, where $u_{ij}$ denotes the  number of elements that are allowed in category $\eta_{ij}$ in the ordinal objective $o^i$. Note that again $0\leq u_{ij}\leq r$ for all $i=1,\dots,p$ and $j\in\{1,\dots,{K_{i}}\}$, and that $\sum_{j\in\{1,\dots,{K_{i}}\}} u_{ij}=r$ for all $i=1,\dots,p$. Therefore, 
it is possible to solve this problem by solving all relevant weighted matroid intersection problems \eqref{p:mi-indep} (with $\M_2$ replaced by $\M_3$) and filtering out all dominated outcome vectors w.r.t.\ the combined ordering relation. In this case,  
the number of 
calls of problems \eqref{p:mi-indep} is bounded by $\bigo( p\cdot r^{\tilde{K}-1})$, where $r$ still denotes the rank of the matroid $\M_1$ and $\tilde{K}=\max\{K_i:i=1,\dots,p\}$. If $p$ and $K_{i}$ are fixed, $i=1,\dots,p$, then the number of scalarized subproblems is polynomially bounded in the input size.
Moreover, in this case every weighted matroid intersection problem can be solved in polynomial time, if we assume that an oracle can determine in polynomial time if a given subset \(I\subseteq E\) is independent or not.

\section{Numerical Results}
\label{sec:numeric}

The Exhaustive Swap Algorithm suggested in \cite{gorski21} as well as the three versions of the Matroid Intersection Algorithm for Ordinal Constraints (Algorithms \ref{alg:mioc}, \ref{alg:mioc-mmpo} and \ref{alg:mioc-mmpc} with the algorithm of Floyd-Warshall) were implemented and numerically tested. As test instances, we consider graphic matroids and partition matroids for $\M_1$, where in the latter case the groundset is partitioned into three subsets. All computations were done on a computer with an Intel(R) Core(TM) i7-7500U CPU 2.70\,GHz processor and 8\,GB RAM. The algorithms were implemented and run in MATLAB, Version R2019b.

In the first experiment we compare the two types of algorithms on a graphic matroid with one sum objective function and one binary objective function. The instances were generated based on random connected undirected graphs $G=(V,E)$ with $n$ nodes and $m$ edges using the implementation of \cite{schn:pkmax:2020}.  
The weight coefficients $w$ of the sum objective are randomly chosen integer values in \(\{1,\ldots, 2\,m\}\), and the values of the binary objective $b$ are random binary values. In both cases we used a uniform distribution. We solved the obtained instances of problem~\eqref{p:bbmp} by the Exhaustive Swap Algorithm and by Algorithm~\ref{alg:mioc-mmpo}.

The numerical results can be found in Table~\ref{t:swapVsMI}. In the first two columns the instance size is given by the number of nodes  and edges $(n,m)$ and the average number of non-dominated outcome vectors $\vert \outNd\vert$ over 100 random instances. The results show clearly that the average running of Algorithm~\ref{alg:mioc-mmpo} increases much faster with the instance size compared to the Exhaustive Swap Algorithm. 

This result is not surprising, because the exhaustive swap algorithm utilizes the specific problem structure, in particular the connectedness of the non-dominated set, as proven in \cite{gorski21}. Nevertheless, the number of problems that are solved with Algorithm~\ref{alg:mioc-mmpo} ($iter$) is quite close to the number of non-dominated outcome vectors ($\outNd$), which indicates that only few redundant problems were solved.

\begin{table}
	\centering\footnotesize
	\setlength{\tabcolsep}{5pt}
	\begin{tabular}{cr|rr|r}
        \toprule
		\multicolumn{2}{c|}{Instanze Size}  & \multicolumn{2}{c|}{ $MIOC^O$} & $ESA$\\
		$(n,m)$&  $\vert \outNd\vert$ &$iter$ & $[s]$ &  $[s]$ \\
		\midrule
		$(7,10)$ & $2.39$& $3.26$& $0.08$& $0.03$\\
		$(7,15)$& $3.32$& $3.84$&  $0.16$& $0.03$\\
		$(7,20)$ & $3.87$ & $3.97$& $0.24$& $0.03$ \\
		$(10,20)$ &$4.36$ &$5.13$  &  $0.54$ & $0.04$\\
		$(10,30)$ & $5.17$ & $5.44$ &  $0.97$& $0.05$  \\
		$(10,40)$ & $5.46$ & $5.53$ & $1.38$ & $0.05$ \\
		$(15,30)$ &$6.05$ & $6.94$& $2.43$ & $0.05$ \\
		$(15,60)$ & $7.77$ & $7.95$ & $6.93$  & $0.07$ \\
		$(15,100)$ & $7.90$ & $7.99$ &  $12.57$ & $0.07$ \\
		$(20,40)$ & $7.58$ & $8.58$ & $7.36$ & $0.07$ \\
		$(20,100)$ & $10.49$ & $10.60$ &  $30.26$ & $0.10$ \\
		$(20,180)$ & $10.35$ & $10.39$ & $59.36$ & $0.11$ \\\bottomrule
	\end{tabular}
	\caption{Average computation time in seconds to solve $100$ instances of problem \eqref{p:bbmp} on a graphic matroid with the Exhaustive Swap Algorithm ($ESA$) and with Algorithm~\ref{alg:mioc-mmpo} ($\mioc^O$).}
	\label{t:swapVsMI}
\end{table}

The strength of all three matroid intersection algorithms for ordinal constraints is that they can be applied to a broader class of problems than the exhaustive swap algorithm, which is restricted to two ordinal categories. In the following tests we used again randomly generated graphs $G=(V,E)$ with $n$ nodes and $m$ edges with objective function coefficients $w$ and $o$. The entries of $w$ and $o$ were generated randomly with uniform distribution in \(\{1,\ldots,2\,m\}\) and  in $\{1,\dots,K\}$ for \(w\) and \(o\), respectively, were $K\in\{3,4,5\}$. The results for $K=3,4,5$ can be found in Tables~\ref{t:MI-K3}, \ref{t:MI-K4} and \ref{t:MI-K5}, respectively. We observe that the number of solutions found for the different problems is quite similar for small problem sizes, but for larger instances and more categories the number of non-dominated points is much smaller for the lexicographic models as compared to the ordinal approach.
The running time depends obviously on the instance size. However, the effect of an increasing number of edges \(m\) is rather limited. A significant influence can be seen by the number of nodes $n$, which determines the rank of the matroid. Furthermore, the number of categories $K$ has an important effect on the running time. 

As expected, reducing the number of considered upper bound vectors $u$ for problem~\eqref{p:mmpo} generally leads to fewer iterations. On average, only little more than half of the iterations are needed in this case. Nevertheless, note that in worst case this strategy may not lead to an improvement. 
In the case of the lexicographic variant~\eqref{p:mmpc} the potential reduction is much more significant. Indeed, the required computation time is drastically reduced in this case, especially for large $K$. For example, for $K=5$ we have a reduction of the running time by a factor of around $20$ in all cases with $n=10$.

We get similar results when testing with partition matroids  rather than graphic matroids. Here, we consider a ground set of $n$ objects and restrict the analysis to partitions of the ground set into three subsets. The upper bounds on the number of elements from each subset are selected such that every basis consists of $\frac{n}{2}$ elements, and the problem is feasible. After defining an instance of a partition matroid $\M_1$ in this way, the objective functions are generated. Each object has an associated weight between $1$ and $10\cdot n$ and is assigned to one of $K$ categories, where  $K\in\{3,4,5\}$. The results for $K=3,4,5$ can be found in Tables~\ref{t:partMI-K3}, \ref{t:partMI-K4} and \ref{t:partMI-K5}, respectively. Again, the improved choice of $u$ leads to significantly better running times. Moreover, the running time increases with the number of elements $n$ and the number of categories $K$.

\begin{table}
	\centering\footnotesize
	\setlength{\tabcolsep}{5pt}
	\begin{tabular}{crrr|rr|rr|rr}
        \toprule
		\multicolumn{4}{c|}{Instance Size} & \multicolumn{2}{c|}{$\mioc$} & \multicolumn{2}{c|}{$\mioc^O$} & \multicolumn{2}{c}{$\mioc^{C_{\min}}$} \\
		$(n,m)$ & $\vert \outNd^O\vert$ &  $\vert \outNd^{C_{\min}}\vert$ &  $\vert \outNd^{C_{\max}}\vert$ & $iter$ & $[s]$  & $iter$ & $[s]$  & $iter$ & $[s]$\\
		\midrule
		$(7,10)$ & $ 3.90$ & $3.65$ & $3.75$  & $28$ & $0.48$ & $17.10$& $0.31$ & $6.30$ & $0.16$ \\
		$(7,15)$ & $7.30$ & $6.05$ & $6.40$ & $28$ & $1.02$ & $18.55$ & $0.77$ & $9.75$ & $0.45$ \\
		$(7,20)$ & $6.20$ & $5.40$ & $5.65$ & $28$ & $1.67$ & $14.45$ & $0.87$ & $8.00$ & $0.56$\\
		$(10,20)$ & $9.10$ & $7.65$ & $7.85$ & $55$ & $5.03$ & $32.50$ & $2.87$ & $12.55$ & $1.47$\\
		$(10,30)$ & $12.55$ & $9.60$ & $11.35$ & $55$ & $9.63$ & $32.80$ & $5.67$ & $13.60$ & $2.61$ \\
		$(10,40)$ & $15.90$ & $11.75$ & $12.95$ & $55$ & $13.56$ & $31.90$ & $8.00$ & $15.25$ & $3.94$\\
		$(15,30)$ & $15.20$ & $11.60$ & $11.40$ & $120$ & $36.16$ & $70.40$ & $21.00$ & $18.60$ & $6.43$\\
		$(15,60)$ & $27.20$ & $18.40$ & $20.65$ & $120$ & $101.30$ & $69.00$ & $58.97$ & $24.30$ & $21.07$\\
		$(15,100)$ & $27.90$ & $18.20$ & $22.65$ & $120$ & $184.26$ & $66.95$ & $105.23$ & $23.10$ & $36.19$\\
		$(20,40)$ & $20.95$ & $14.65$ & $15.40$ & $210$ & $158.48$ & $113.85$ & $85.85$ & $23.95$ & $20.53$\\
		$(20,100)$ & $38.45$ & $24.55$ & $27.00$ & $210$ & $578.96$ & $101.00$ & $282.43$ & $30.95$ & $87.42$\\
		$(20,180)$ & $46.55$ & $27.20$ & $33.45$ & $210$ & $1\,162.27$ & $115.35$ & $648.08$ & $33.70$ & $190.49$\\\bottomrule
	\end{tabular}
	\caption{Numerical results for a graphic matroid and $K=3$ categories. For every problem size $20$ instances were solved to obtain average results.}
	\label{t:MI-K3}
\end{table}

\begin{table}
	\centering\footnotesize
	\setlength{\tabcolsep}{5pt}
	\begin{tabular}{crrr|rr|rr|rr}
        \toprule
		\multicolumn{4}{c|}{Instance Size} & \multicolumn{2}{c|}{$\mioc$} & \multicolumn{2}{c|}{$\mioc^O$} & \multicolumn{2}{c}{$\mioc^{C_{\min}}$} \\
		$(n,m)$ & $\vert \outNd^O\vert$ &  $\vert \outNd^{C_{\min}}\vert$ &  $\vert \outNd^{C_{\max}}\vert$ & $iter$ & $[s]$  & $iter$ & $[s]$  & $iter$ & $[s]$\\
		\midrule
		$(7,10)$ & $4.15$ & $3.70$ & $3.70$  & $84$ & $1.34$ & $44.45$ & $0.73$ & $9.00$ & $0.25$ \\
		$(7,15)$ & $7.50$ & $5.85$ & $6.25$ & $84$ & $2.83$ & $43.85$ & $1.61$ & $12.10$ & $0.56$\\
		$(7,20)$ & $10.50$ & $7.85$ & $8.10$ & $84$ & $4.44$ & $52.05$ & $2.88$ & $14.80$ & $1.93$ \\
		$(10,20)$ & $12.15$ & $8.65$ & $8.60$ & $220$ & $18.15$ & $120.25$ & $9.38$ & $18.75$ & $2.01$\\
		$(10,30)$ & $20.60$ & $12.40$ & $14.80$ & $220$ & $36.78$ & $132.95$ & $22.01$ & $23.85$ & $4.46$\\
		$(10,40)$ & $23.70$ & $14.35$ & $17.35$ & $220$ & $53.79$ & $118.00$ & $28.27$ & $23.90$ & $5.95$\\
		$(15,30)$ & $29.05$ & $16.80$ & $17.55$ & $680$ & $191.27$ & $395.75$ & $107.14$ & $37.15$ & $12.84$\\
		$(15,60)$ & $56.00$ & $24.95$ & $34.10$ & $680$ & $560.72$ & $350.15$ & $291.72$ & $42.40$ & $35.97$\\
		$(15,100)$ & $63.60$ & $28.00$ & $38.70$ & $680$ & $1\,042.33$ & $377.45$ & $597.78$ & $44.15$ & $69.33$ \\
		$(20,40)$ & $44.80$ & $21.80$ & $23.50$ & $1\,540$ & $1\,073.17$ & $708.20$ & $502.36$& $43.80$ & $37.32$\\\bottomrule
	\end{tabular}
	\caption{Numerical results for a graphic matroid and $K=4$ categories. For every problem size $20$ instances were solved to obtain average results.}
	\label{t:MI-K4}
\end{table}

\begin{table}
	\centering\footnotesize
	\setlength{\tabcolsep}{5pt}
	\begin{tabular}{crrr|rr|rr|rr}
        \toprule
		\multicolumn{4}{c|}{Instance Size} & \multicolumn{2}{c|}{$\mioc$} & \multicolumn{2}{c|}{$\mioc^O$} & \multicolumn{2}{c}{$\mioc^{C_{\min}}$} \\
		$(n,m)$ & $\vert \outNd^O\vert$ &  $\vert \outNd^{C_{\min}}\vert$ &  $\vert \outNd^{C_{\max}}\vert$ & $iter$ & $[s]$  & $iter$ & $[s]$  & $iter$ & $[s]$\\
		\midrule
		$(7,10)$ & $ 4.20$ & $3.85$ & $3.85$  & $210$ & $3.10$ & $106.75$& $1.57$ & $10.95$ & $0.27$ \\
		$(7,15)$ & $10.75$ & $7.55$ & $8.10$ & $210$ & $6.65$ & $121.75$ & $3.88$ & $18.05$ & $0.81$\\
		$(7,20)$ & $14.90$ & $10.00$ & $11.25$ & $210$ & $10.56$ & $124.30$ & $6.14$ & $20.35$ & $1.36$\\
		$(10,20)$ & $22.85$ & $11.75$ & $12.65$ & $715$ & $57.29$ & $429.80$ & $33.25$ & $30.30$ & $3.19$ \\
		$(10,30)$ & $24.70$ & $14.55$ & $15.95$ & $715$ & $114.37$ & $396.00$ & $62.16$ & $32.00$ & $5.57$\\
		$(10,40)$ & $30.10$ & $15.65$ & $19.80$ & $715$ & $175.93$ & $410.75$ & $98.27$ & $32.80$ & $8.22$ \\\bottomrule
	\end{tabular}
	\caption{Numerical results for a graphic matroid and $K=5$ categories. For every problem size $20$ instances were solved to obtain average results.}
	\label{t:MI-K5}
\end{table}

\begin{table}
	\centering\footnotesize
	\setlength{\tabcolsep}{5pt}
	\begin{tabular}{crrr|rr|rr|rr}
		\toprule
		\multicolumn{4}{c|}{Instance Size} & \multicolumn{2}{c|}{$\mioc$} & \multicolumn{2}{c|}{$\mioc^O$} & \multicolumn{2}{c}{$\mioc^{C_{\min}}$} \\
		$n$ & $\vert \outNd^O\vert$ &  $\vert \outNd^{C_{\min}}\vert$ &  $\vert \outNd^{C_{\max}}\vert$ & $iter$ & $[s]$  & $iter$ & $[s]$  & $iter$ & $[s]$\\
		\midrule
		$10$ & $2.55$ & $2.50$ & $2.50$  & $21$ & $0.23$ & $11.20$& $ 0.10$& $4.60$ & $0.08$ \\
		$20$ & $4.75$ & $4.40$ & $4.45$  & $66$ & $3.26$ & $35.45$ & $1.83$ & $8.40$ & $0.63$ \\
		$30$ & $9.90$ & $8.00$ & $8.40$  & $136$ & $20.26$ & $76.50$ & $11.43$ & $13.45$ & $2.68$ \\
		$40$ & $16.85$ & $12.55$ & $12.75$  & $231$ & $81.26$ & $127.25$ & $45.14$ & $20.40$ & $8.92$ \\
		$50$ & $24.05$ & $17.60$ & $15.90$  & $351$ & $229.91$ & $201.90$ & $131.29$ & $26.30$ & $21.43$ \\
		$60$ & $34.00$ & $21.10$ & $21.75$  & $496$ & $552.54$ & $284.20$ & $318.62$ & $32.80$ & $46.41$ \\
		$70$ & $39.70$ & $24.50$ & $25.40$  & $666$ & $1\,177.52$ & $375.90$ & $662.34$ & $37.30$ & $83.06$ \\\bottomrule
	\end{tabular}
	\caption{Numerical results for a partition matroid with three subsets and $K=3$ categories. For every problem size $20$ instances were solved to obtain average results.}
	\label{t:partMI-K3}
\end{table}

\begin{table}
	\centering\footnotesize
	\setlength{\tabcolsep}{5pt}
	\begin{tabular}{crrr|rr|rr|rr}
		\toprule
		\multicolumn{4}{c|}{Instance Size} & \multicolumn{2}{c|}{$\mioc$} & \multicolumn{2}{c|}{$\mioc^O$} & \multicolumn{2}{c}{$\mioc^{C_{\min}}$} \\
		$n$ & $\vert \outNd^O\vert$ &  $\vert \outNd^{C_{\min}}\vert$ &  $\vert \outNd^{C_{\max}}\vert$ & $iter$ & $[s]$  & $iter$ & $[s]$  & $iter$ & $[s]$\\
		\midrule
		$10$ & $2.55$ & $2.45$ & $2.40$  & $56$ & $0.41$ & $31.10$ & $0.23$ & $6.15$ & $0.10$ \\
		$20$ & $7.10$ & $5.80$ & $5.85$  & $286$ & $12.60$ & $173.90$ & $7.43$ & $14.05$ & $0.94$ \\
		$30$ & $15.15$ & $9.75$ & $10.05$  & $816$ & $117.69$ & $451.70$ & $62.22$ & $22.85$ & $4.33$ \\
		$40$ & $40.40$ & $20.95$ & $21.25$  & $1\,771$ & $594.17$ & $967.80$ & $317.36$ & $42.25$ & $18.51$ \\
		$50$ & $44.90$ & $22.50$ & $22.20$  & $3\,276$ & $2\,081.56$ & $1\,946.25$ & $1\,191.88$ & $49.70$ & $41.17$ \\\bottomrule
	\end{tabular}
	\caption{Numerical results for a partition matroid with three subsets and $K=4$ categories. For every problem size $20$ instances were solved to obtain average results.}
	\label{t:partMI-K4}
\end{table}

\begin{table}
	\centering\footnotesize
	\setlength{\tabcolsep}{5pt}
	\begin{tabular}{crrr|rr|rr|rr}
		\toprule
		\multicolumn{4}{c|}{Instance Size} & \multicolumn{2}{c|}{$\mioc$} & \multicolumn{2}{c|}{$\mioc^O$} & \multicolumn{2}{c}{$\mioc^{C_{\min}}$} \\
		$n$ & $\vert \outNd^O\vert$ &  $\vert \outNd^{C_{\min}}\vert$ &  $\vert \outNd^{C_{\max}}\vert$ & $iter$ & $[s]$  & $iter$ & $[s]$  & $iter$ & $[s]$\\
		\midrule
		$10$ & $3.65$ & $3.15$ & $3.25$  & $126$ & $0.84$ & $71.55$ & $0.47$ & $8.20$ & $0.14$ \\
		$20$ & $10.65$ & $8.20$ & $7.75$  & $1\,001$ & $43.34$ & $539.15$ & $22.62$ & $23.20$ &$1.67$  \\
		$30$ & $28.25$ & $15.35$ & $16.35$  & $3\,876$ & $534.59$ & $2\,293.50$ & $308.62$ & $41.50$ & $7.83$ \\\bottomrule
	\end{tabular}
	\caption{Numerical results for a partition matroid with three subsets and $K=5$ categories. For every problem size $20$ instances were solved to obtain average results.}
	\label{t:partMI-K5}
\end{table}

Note that further speed-ups can be expected by parallel implementations of the matroid intersection algorithms.

\section{Conclusion}
\label{sec:concl}
In this paper we consider single- and multi-objective matroid optimization problems that combine ``classical'' sum objective functions with one or several ordinal objective functions. Besides the concept of ordinal optimality, we consider two variants of lexicographic optimization that lexicographically maximize the number of ``good'' elements or minimize the number of ``bad'' elements, respectively. In the case of (single-objective) ordinal optimization, we show that these concepts are actually equivalent for matroids, and that optimal solutions can be found by a simple and efficient greedy strategy. In the multi-objective setting, we use variants of $\eps$-constraint scalarizations to obtain a polynomial number of matroid intersection problems, from which the non-dominated sets of the respective problems can be derived by simple filtering operations. This yields an overall polynomial-time algorithm for multi-objective ordinal matroid optimization problems. Numerical tests on graphic matroids and on partition matroids validate the efficiency of this approach.

Future research should focus on a further analysis of the similarities and differences between multi-objective optimization problems with classical sum objectives and with ordinal objectives. Moreover, alternative (partial) orderings may be considered and analyzed in the light of different scalarization techniques.

\section*{Acknowledgment}
The authors thankfully acknowledge financial support by the Deutsche Forschungsgemeinschaft, project number~KL~1076/11-1.

\bibliography{literaturMatroid}
\end{document}